\documentclass[11 pt]{article}
\usepackage{amssymb,amsmath,amsfonts,amsthm}

\providecommand{\U}[1]{\protect\rule{.1in}{.1in}}
\newtheorem{theorem}{Theorem}[section]
\newtheorem{lemma}{Lemma}[section]

\newtheorem{Rem}{Remark}[section]
\theoremstyle{remark}

\usepackage{color}
\def\1{{{\mbox{${\rm{1\negthinspace\negthinspace I}}$}}}}

\newcommand{\eref}[1]{(\ref{#1})}

\newcommand\beq{\begin{equation}}
\newcommand\eeq{\end{equation}}

\begin{document}

\title{Rates of convergence  in the central limit theorem for the elephant random walk with random step sizes }

\author{J\'er\^ome Dedecker\footnote{Universit\'{e} Paris Cit\'e, Laboratoire MAP5 and CNRS UMR 8145, 75016
Paris, France.}\footnotemark[1],\ \   Xiequan Fan\footnote{Center for Applied Mathematics,
Tianjin University, Tianjin 300072, China.}\footnotemark[2],
\ \ Haijuan Hu\footnote{School of Mathematics and Statistics, Northeastern University at Qinhuangdao, Qinhuangdao, China.}\footnotemark[3],
\ \ Florence Merlev\`{e}de\footnote{LAMA, Univ Gustave Eiffel, Univ Paris Est Cr\'{e}teil, UMR 8050 CNRS, F-77454 Marne-La-Vall\'{e}e, France.}\footnotemark[4]
}

\date{}

\maketitle

\abstract{In this paper, we consider a generalization of the elephant random walk  model. Compared to the usual elephant random walk,
an interesting feature of this model is that the step sizes form a sequence of positive  independent and identically distributed random variables   instead of a fixed constant. For this model, we establish the law of the iterated logarithm, the central limit theorem,  and  we obtain rates of convergence in the central limit theorem with  respect to the  Kologmorov, Zolotarev and Wasserstein distances. We emphasize that,  even in case of the usual   elephant random walk, our results concerning the rates of convergence in the central limit theorem are new. }

\medskip

\noindent {\bf Keywords.} Elephant random walk, law of  the iterated logarithm, normal approximations, Wassertein's  distance,
 Berry-Esseen bounds, central limit theorem
\medskip

\noindent {\bf Mathematics Subject Classification (2010):} 60G42, 60F05, 60E15, 82C41

\section{Introduction}

 The elephant random walk (ERW) is a type of one-dimensional random walk on integers, which has a complete memory of its whole history.
   The  ERW was introduced in 2004  by Sch{\"u}tz and Trimper \cite{schutz2004elephants}  in order to study the memory effects of a non-Markovian random walk. The model has a link to a famous saying that elephants can always remember where they have been.

    The   ERW can be defined as  follows.
It starts at time $n=0$, with  position $T_0=0$. At time $n=1$, the elephant
 moves to $1$ with probability $q$ and to $-1$ with probability $1-q$, where $q \in (0, 1].$
 So the position of the elephant at time $n=1$ is given by $T_1=X_1,$ with $X_1$
 a Rademacher $\mathcal{R}(q)$ random variable.
 At time $n+1,$   $n\geq 1,$ the step $X_{n+1}$ is determined stochastically by the following rule.
 Let $n'$ be an integer  which is  chosen from the set $\{1, 2,\ldots, n\}$
 uniformly at random.
 If $X_{n'}=1$, then
 \begin{displaymath}
 X_{n+1} =   \left\{ \begin{array}{ll}
1 & \textrm{with  probability  $p$}\\
-1 & \textrm{with  probability  $1-p$}.
\end{array} \right.
\end{displaymath}
 If $X_{n'}=-1$, then
 \begin{displaymath}
 X_{n+1} =   \left\{ \begin{array}{ll}
1 & \textrm{with  probability  $1-p$}\\
-1 & \textrm{with  probability  $p$}.
\end{array} \right.
\end{displaymath}
Equivalently,  $X_{n+1}$ is determined stochastically by the following rule:
 \begin{displaymath}
 X_{n+1 } =   \left\{ \begin{array}{ll}
X_{n'} & \textrm{with  probability  $p$}\\
-X_{n'} & \textrm{with  probability  $1-p$}.
\end{array} \right.
\end{displaymath}
The last equality suggests that at time $n+1,$ we reinforce $X_{n'}$ with  probability $p$ and reduce $X_{n'}$ with  probability $1-p$.
Thus, for $n\geq 2,$  the position of the elephant at time $n $ is
 \begin{eqnarray}\label{adfs0}
  T_{n }=\sum_{i=1}^{n }X_i, \text{ where }
 X_{n }=\alpha_{n}X_{\beta_{n}},
\end{eqnarray}
with $\alpha_n$ has a Rademacher distribution $\mathcal{R}(p)$, $p\in (0, 1],$ and $\beta_n$ is uniformly distributed over the integers $\{1, 2,\ldots, n-1\}$. Moreover, $\alpha_n$ is independent of $(X_i)_{1\leq i \leq n}$, and the random variables  $(\beta_i)_{i\geq1}$  are independent.
Here  $p $ is called  the memory parameter. The ERW is respectively  called diffusive,  critical and   superdiffusive according to  $p \in (0, 3/4),$ $p=3/4$ and $p \in (3/4, 1]$.

 %

 The description of the asymptotic behavior of the ERW has
motivated many interesting works.
By showing the connection of the ERW model with Polya-type urns, Baur and Bertoin obtained the functional form of the  central  limit theorem.
Coletti, Gava and Sch{\"u}tz \cite{coletti2017central,CGS17} proved the central limit theorem (CLT)  and a  strong invariance principle for $p \in (0,  3/4]$ and a law of large numbers for $p \in (0, 1).$ Moreover, they also showed that if $p \in (3/4, 1]$, then the ERW converges to a non-degenerate random variable which is not normal.
V\'{a}zquez Guevara \cite{G19} gave the almost sure CLT.  Bercu \cite{B18} recovered the CLT via a martingale method. Bertoin \cite{B21}
  studied how memory impacts passages at the origin for the ERW.
Recently,   a number of variations for the ERW has been introduced. For instance, Bercu and Lucile \cite{BL19} introduced a multi-dimensional type ERW, and  established the multivariate CLT. Gut and Stadtm\"{u}ller \cite{GS71} investigated the case when the elephant remembers only a finite part of the first or last steps. 
 Laulin \cite{L22} studied the asymptotic analysis for the reinforced ERW model. Recently, Bercu \cite{B22}  considered the ERW with stops playing hide and seek with the Mittag-Leffler distribution.

In this paper, we consider a generalization of ERW, called ERW with random step sizes, such that the step sizes are random and varying in time.
The ERW with random step sizes was introduced  by  Fan  and   Shao \cite{FS23}.
 Let $(Z_i)_{i\geq1}$ be a sequence of  positive  independent and identically distribution  (i.i.d.) random variables, with finite mean  $\nu=  \mathbf{E} Z_1$ and variance $\textrm{Var}(Z_1)=\sigma^2\geq 0$. Moreover, $(Z_i)_{i\geq1}$ is independent of $(X_i)_{i\geq1}$. An ERW with  random step sizes may be described as follows. At time $n = 1$, the elephant
 moves to $Z_1$ with probability $q$ and to $-Z_1$ with probability $1-q$.
 So the position $Y_1$ of the elephant at time $n=1$ is given by   the following rule:
 \begin{displaymath}
 Y_{1} =   \left\{ \begin{array}{ll}
Z_1 & \textrm{\ \ \ with  probability  $q$}\\
-Z_1 & \textrm{\ \ \ with  probability  $1-q$}.
\end{array} \right.
\end{displaymath}
For $n\geq2,$  instead of (\ref{adfs0}), the position of the elephant with random step sizes at time $n$ is
\[S_{n}=\sum_{i=1}^{n}Y_i, \text{ where }
  Y_{n}=\alpha_{n  }X_{\beta_{n  }}Z_{n}.\]
Notice that $|Y_{n}|=Z_{n}$ for all $n\geq 1$. Thus at time $n$, the step size is $Z_{n} $ which is a random variable.
Without loss of generality, we may assume that $\nu=1$ (otherwise, we may consider the case $S_{n}/\nu$ instead of $S_{n}$).
Clearly, when $\sigma^2=0,$ we have $Z_{n}\equiv1$ a.s.\ and then the ERW with  random step sizes reduces to the usual ERW.

In this paper,  we are interested  in establishing the  law of the  iterated logarithm and  convergence rates for normal approximations of $S_{n}$ in terms of Kolmogorov distances and  Zolotarev  or  Wasserstein distances. As we shall see,  even in case of the usual  ERW, our results concerning the rates of convergence in the central limit theorem are new.

Throughout the paper,  $C$  will denote a finite and positive constant that is allowed to depend on some fixed parameters such as $p$, $q$, $r$, $\rho$, $\gamma$, ${\mathbf E} Z_1^{2+ \rho}$ and so on,  but not on $n$. This constant may vary from line to line.    For  two sequences of positive numbers $(a_n)_{n\geq 1}$ and $(b_n)_{n\geq 1}$,   write $a_n \asymp b_n$ if there exists a constant $C>0$ such that ${a_n}/{C}\leq b_n\leq C a_n$ for all sufficiently large $n$. We shall also use the notation $a_n \ll b_n$ to mean that there exists a positive constant $C$ not depending on $n$ such that $a_n \leq C b_n$.
We also  write $a_n \sim b_n$ if  $\lim_{n\rightarrow \infty} a_n/b_n =1$.
We denote by $\mathcal{N}(0, \sigma^2)$   the normal  distribution with mean $0$ and variance $\sigma^2$.
$\mathcal{N}$  will designate a standard normal random variable, and we will denote by $\Phi(\cdot)$  the cumulative distribution function of a standard normal random variable. We shall use the notations
$\varphi$ to denote the density of a $\mathcal{N}(0, 1)$ distribution and $\varphi_{\sigma^2}$ to denote the density of a $\mathcal{N}(0, \sigma^2)$ distribution.

\section{Law of iterated logarithm}
\setcounter{equation}{0}
The almost sure convergence  for the ERW has been established by  Bercu \cite{B18} and  Bercu and Lucile \cite{BL19} via a martingale method.
In the next theorem, we give the corresponding results for the ERW with random step sizes.

\begin{theorem}\label{th3.1}  Assume that $p \in (0,  1] $  and  $ {{\mathbf E}} Z_1^2 < \infty.$
 \begin{description}
  \item[\textbf{[i]}]  If $p \in (0,  3/4)$,  then
\begin{equation}\label{fsv12}
\limsup_{n \rightarrow \infty} \frac{|S_{n}|}{\sqrt{n \log \log n} } \leq  \sqrt{\frac{2}{3-4p} }+ \sqrt{2} \sigma\ \ \ \ \textrm{a.s.}
\end{equation}

  \item[\textbf{[ii]}]  If $p = 3/4$,  then
\begin{equation}\label{fsv13}
\limsup_{n \rightarrow \infty} \frac{|S_{n}|}{\sqrt{n \log n  \log \log \log n } } = \sqrt{2}  \ \ \ \ \textrm{a.s.}
\end{equation}

\item[\textbf{[iii]}]  If $p \in ( 3/4,   1 ]$,  then
\begin{equation}\label{fsv22}
\lim_{n \rightarrow \infty} \frac{S_{n}}{ n^{2p-1} } = L \ \ \ \ \textrm{a.s.},
\end{equation}
where $L$ is a non-degenerate and non Gaussian random variable.
\end{description}
\end{theorem}
%
%

We immediately deduce from Theorem \ref{th3.1} that the  law of large numbers  holds.
 If $p \in (0,  1)$,  then
\begin{equation}\nonumber
\lim_{n \rightarrow \infty} \frac{S_{n}}{n } = 0\ \ \ \ \textrm{a.s.}
\end{equation}
Hence, we find again the Kolmogorov strong  law of large numbers given by Coletti,  Gava and Sch{\"u}tz \cite{coletti2017central}. See also   Bercu \cite{B18} and Bercu  and  Lucile \cite{BL19} for closely  related results. In particular,
 Bercu \cite{B22} gave the law of the  iterated logarithm for ERW
 with stops playing hide and seek with the Mittag-Leffler distribution.

 \begin{Rem}
 In \cite{CGS17} a strong invariance principle is proved for the usual ERW, from which the compact law of the iterated logarithm can be easily deduced. We would like to mention that the rate of approximation in this strong invariance principle can be improved by using Theorem 2.1 of Shao \cite{shao93}. More precisely, letting $v_n= \sum_{i=1}^n a_i^2$ and applying Shao's results, one can prove that
 \begin{description}
  \item[\textbf{[i]}]  If $p\leq 5/8$,    there exists a standard Wiener process $W{(t)}$ such that, for any $\delta >0$,
$$
    a_nT_n -W(v_n) = o\left(a_n n^{1/4} (\log n)^{3/4} (\log \log n)^{\delta + 1/4}\right) \  \text{a.s.}
$$

  \item[\textbf{[ii]}]  If $5/8<p < 3/4$,  then there exists a standard Wiener process $W{(t)}$ such that
$$
    a_nT_n -W(v_n) = o\left(u_n  (\log n)^{1/2} \right) \  \text{a.s., where $u_n$ is any sequence increasing to $\infty$.}
$$

\item[\textbf{[iii]}]  If $p=3/4$,  then there exists a standard Wiener process $W{(t)}$ such that
$$
    a_nT_n -W(v_n) = o\left(u_n  (\log \log n)^{1/2} \right) \  \text{a.s., where $u_n$ is any sequence increasing to $\infty$.}
$$
\end{description}
 \end{Rem}

\section{Normal approximations }\label{secm}
\setcounter{equation}{0}
In this section, we consider the normal approximation for $S_n.$
For $p \in (0, 3/4] $ and $n\geq 2,$ denote
$$ a_1=1, \ \ \ \ \ \  a_n=\frac{\Gamma(n)\Gamma(2p)}{\Gamma(n+2p-1)}  \ \ \ \ \ \ \textrm{and} \ \ \ \ \ \   v_n=\sum_{i=1}^n a_i^2.$$
Here $\Gamma(s)=\int_0^{\infty}t^{s-1}e^{-t}dt, s>0, $ is the Gamma function.
Notice that  $a_n$ and $v_n$ are both positive. Moreover, by taking into account the relations (\ref{a15}) and (\ref{a16}) of Section \ref{secdp},
 it is known that
$v_n \rightarrow \infty$ as $n\rightarrow \infty $.
\subsection{Central limit theorem}

\begin{theorem} \label{th0}  Assume that  $ {\mathbf E} Z_1^2 < \infty   $.
Then,  for   $p \in (0, 3/4],$
\begin{equation}\label{CLT}
\displaystyle   \frac{{\tilde S}_n}{\displaystyle  \sqrt{v_n+ n a_n^2 \sigma^2}}  \stackrel{\mathbf{D}}{\longrightarrow }   \mathcal{N},
\end{equation}
where $ {\tilde S}_n=a_nS_n -(2q-1) $ and $\stackrel{\mathbf{D}}{\longrightarrow}$ stands for convergence in distribution.
\end{theorem}

When $p \in (0, 3/4]$ and $\sigma^2=0,$ the central limit theorem for the usual ERW  has been established by Coletti, Gava and Sch{\"u}tz \cite{coletti2017central,CGS17}  and Bercu \cite{B18}. Theorem  \ref{th0} then extends the CLT to the case $\sigma^2 >0$.

Note also that since $v_n \rightarrow \infty$ as $n\rightarrow \infty $, \eqref{CLT} also holds if we replace  $ {\tilde S}_n$ by $a_n S_n$ (indeed the term
$(2q-1) /  \sqrt{v_n+ n a_n^2 \sigma^2}$ goes to $0$ as $n \rightarrow \infty$).  However, it is convenient to consider random variables which are centered as
 $ {\tilde S}_n$. Indeed if we consider $a_n S_n$ rather than $ {\tilde S}_n$ then the corresponding rates of convergence in the CLT are slower (see for instance  Remark \ref{KDNC}  concerning the Kolmogorov distance).

 Concerning the normalizing sequence, note that when $p \in (0,3/4)$, then (\ref{a10}) combined with (\ref{a15}) entails that
 \[
b_n :=  a_n^{-1} \sqrt{v_n+ n a_n^2 \sigma^2} \sim \sqrt{n (\sigma^2 + 1/(3-4p)) },
 \]
whereas when $p=3/4$,   (\ref{a10}) combined with (\ref{a15})  implies that $b_n \sim \sqrt{  n \log n  }$.

Finally, notice that when $p \in (3/4, 1]$, $v_n+ n a_n^2 \sigma^2$  is bounded and $a_n$ is in order of $n^{1-2p}$ (see \eqref{a10}).
  Thus, by  \eref{fsv22},    $ {\tilde S}_n /  \sqrt{v_n+ n a_n^2 \sigma^2} $   converges to a non-degenerate and non Gaussian random  random variable, and
  \eref{CLT} does not hold. Hence, we  shall only   consider the case $p \in
 (0, 3/4]$ in the remaining of the paper.

\subsection{Convergence rates in terms of Kologmorov's distance}
For a real-valued random variable $X,$  denote
$$\mathbf{K}(X )=  \sup_{u \in \mathbf{R}}\Big| \mathbf{P} \big(X  \leq u \big)  -  \Phi (u)\Big|$$
the Kologmorov's distance between the law of $X$ and the standard normal distribution.

We start by giving the rates of convergence in the central limit theorem for the \mbox{Kolmogorov} distance and the usual ERW.
\begin{theorem} \label{thmBEERW}
Let $p \in (0, 3/4].$
Then 
 the following inequality holds:
\[
\mathbf{K}\Bigg(  \frac{a_nT_n -(2q-1)}{ \sqrt{v_n}} \Bigg)  \ll \frac{1}{\sqrt{n}} + \frac{1}{v_n} .
\]
\end{theorem}
Taking into account \eqref{a15} and \eqref{a16} of the next section, it follows that when $p \leq 5/8$, the rate in Theorem \ref{thmBEERW} is $1/{\sqrt n}$ and when $p >5/8$, the rate
is $1/v_n$.
Note also that Theorem \ref{thmBEERW}  improves on the previous results obtained by Fan \emph{et al.}\ \cite{FHM20} who obtained the rate
$(\log n)/\sqrt{n}$ when $p \leq 1/2$,  $(\log n)/ \sqrt{v_n}$ when $p \in (1/2, 3/4)$ and $(\log \log  n)/ \sqrt{v_n}$ when $p= 3/4$.

\begin{Rem} \label{KDNC}
Using Lemma \ref{LEMMA-APX-1}, we then derive that for the non centered quantity $a_nT_n/ \sqrt{v_n}$ we have  the following rate of convergence in the CLT
\[
\mathbf{K}\Bigg(  \frac{a_nT_n }{ \sqrt{v_n}} \Bigg)  \ll \frac{1}{\sqrt{n}} + \frac{1}{\sqrt{v_n}} .
\]
\end{Rem}

The following theorem gives a Berry-Esseen type bound for the ERW with random step sizes.
Theorem  \ref{thmBEERW} can be viewed as a consequence of it (taking $\rho=1$ and $\sigma=0$.)

\begin{theorem} \label{thm6}
Assume that $M(\rho)=\mathbf{E} Z_1^{2+\rho}< \infty$ for a constant $ \rho \in (0, 1]. $
The following inequality holds:
\[
\mathbf{K}\Bigg(  \frac{a_nS_n -(2q-1)}{ \sqrt{v_n+ n a_n^2 \sigma^2}} \Bigg)  \ll  \frac{1}{ n ^{ \rho/2 } } + \frac{1}{v_n} .
\]
\end{theorem}


\subsection{Convergence rates in terms of Zolotarev's distances}
Zolotarev's distance of order $r$ can be described as follows: for two probability measures $\mu$ and $\nu,$ and any $r>0,$
\begin{eqnarray*}
	\zeta_r(\mu, \nu)= \sup  \bigg \{   \int f d \mu - \int f d \nu: \    f \in \Lambda_r      \bigg\},
\end{eqnarray*}
where $\Lambda_r $ is defined as follows: denoting by $l$ the natural integer satisfying $l < r \leq l+1$, $\Lambda_r$
is the class of real functions $f$ which are $l$-times continuous differentiable and such that
$$   |f^{(l)}(x) - f^{(l)}(y) |   \leq |x-y|^{r-l}\ \ \ \ \textrm{for any }(x, y) \in \mathbf{R}^2.$$
For two random variables $X$ and $Y$, let $\mu_X$ and  $\nu_Y$ be respectively the probability laws of  $X$ and $Y$. We  denote
 $\zeta_r (X, Y) = \zeta_r (\mu_X, \nu_Y)$ to soothe the notations.
 \medskip

We start by giving the rates of convergence in the central limit theorem in terms of  Zolotarev distances for  the usual ERW.

\begin{theorem}\label{thm1ERW}
 Let $r \in (0, 2].$
 Then 
 the following inequality holds:
\[
\zeta_r \bigg( \frac{a_nT_n -(2q-1)}{ \sqrt{v_n}} , \  \mathcal{N}  \bigg)
\ll  \frac{1}{n^{r /2}}   +   \frac{1}{v_n } .
\]
  \end{theorem}
  Note that $\zeta_1=W_1$, where $W_1$ is the Wasserstein distance of order 1 (see the next section for the definition). Hence, for $r=1$,
  Theorem \ref{thm1ERW} improves on the previous results obtained by Ma \textit{et al.}\ \cite{MEF22}.

 \medskip

The following theorem gives some  bounds for the ERW  with random step sizes.

\begin{theorem}\label{thm1}
 Assume that $\mathbf{E} Z_1^{2+\rho}< \infty$ for a constant $ \rho \in (0, 1]. $ Let $r \in (0, 2].$
The following inequality holds:
\[
 \zeta_r \bigg( \frac{a_nS_n -(2q-1)}{ \sqrt{v_n+ n a_n^2 \sigma^2}} , \  \mathcal{N}  \bigg)
\ll  \frac{1}{n^{(r\wedge \rho) /2}}   +   \frac{1}{v_n } .
\]
\end{theorem}

\subsection{Convergence rates in terms of Wasserstein's distances}

We first recall the definitions of Wasserstein's distances.  Let $\mathcal{L}(\mu, \nu)$ be a set
of probability laws on $\mathbf{R}^2$ with marginals $\mu$ and  $\nu$.  The  Wassertein distance of order $r >0 $ is defined
as follows:
\begin{eqnarray*}
	W_r(\mu, \nu)= \inf  \bigg \{ \Big(\int |x-y|^r \mathbf{P}(dx, dy) \Big)^{1/ \max(1, r)} :\ \mathbf{P} \in \mathcal{L}(\mu, \nu)      \bigg\} .
\end{eqnarray*}
For two random variables $X$ and $Y$ with respective laws $\mu_X$ and  $\nu_Y$ ,  denote
$W_r(X, Y) =W_r(\mu_X, \nu_Y) $ to soothe the notation.

For $r \in (0, 1]$, using Kantorovich-Rubinstein's theorem with the metric $d(x, y)=|x-y|^r,$ we have
$
W_r(\mu, \nu) = \zeta_r(\mu, \nu)$. Hence, for $r \in (0, 1]$, the inequality of Theorem \ref{thm1} holds for the distance $W_r$.

\smallskip

Next, for $r \in [1,2]$, from  \cite{Rio2009}, we know that there exists a positive constant $c_r$ such that
\begin{equation} \label{Rioine}
W_r(\mu, \nu) \leq  c_r (\zeta_r(\mu, \nu))^{1/r}.
\end{equation}
Therefore, starting from Theorem \ref{thm1ERW}, we derive for the usual ERW that for any $r \in [1,2]$,
 \[
W_r \bigg( \frac{a_nT_n -(2q-1)}{ \sqrt{v_n}} , \  \mathcal{N}  \bigg)
\ll  \frac{1}{ \sqrt{n}}   +   \frac{1}{v^{1/r}_n } .
\]
Combining this last estimate with some precise upper bounds for independent random variables given by Bobkov \cite{Bob18}, we get the following result for the ERW with random step sizes.

\begin{theorem} \label{thm2}
 Let $r \in [1,2]$, and assume that
 $ \mathbf{E} Z_1^{2+\rho}< \infty$ for a constant $ \rho \in (0, r]$. The following inequalities hold:
\[
W_r \bigg( \frac{a_nS_n -(2q-1)}{ \sqrt{v_n}} , \  \mathcal{N}  \bigg)
\ll  \frac{1}{n^{\rho /2r}}   +   \frac{1}{v_n^{1/r} } .
\]
\end{theorem}
Note that when $\rho \in (0,1]$, Theorem \ref{thm2} is a direct consequence of Theorem \ref{thm1} combined with  inequality \eqref{Rioine}.

\section{Preliminary considerations and lemmas} \label{secdp}
\setcounter{equation}{0}

By the well-known Stirling  formula $$\log  \Gamma(x)  =\left (x-\frac12 \right ) \log x -x + \frac12 \log 2 \pi + O\left (\frac{1}{x} \right )  \ \ \ \  \ \textrm{as}\ \ x \rightarrow \infty,$$ we deduce that for $p \in (0, 1],$
\begin{eqnarray}\label{a10}
\lim\limits_{n\to \infty} a_n n^{ 2p-1} =\Gamma(2p).
\end{eqnarray}
Moreover, for $  p\in (0, 3/4)$, we have
\begin{eqnarray}\label{a15}
\lim\limits_{n\to \infty} \frac{ v_n}{ n^{3-4p}} = \frac{\Gamma{(2p)}^2}{3-4p},
\end{eqnarray}
and, for $ p=3/4$,
\begin{eqnarray}\label{a16}
\lim\limits_{n\to \infty} \frac{ v_n}{  \log n}    =\frac{\pi}{4}.
\end{eqnarray}
See also Bercu \cite{B18} for the equalities (\ref{a10})-(\ref{a16}).
Denote $$\gamma _n=1+\frac{2p-1}{n},\ \ \ \ \ \ \ \ \ n\geq 1.  $$
It is easy to see that $a_1=1$ and for $n\geq 2,$
$$a_n= \prod_{i=1}^{n-1}\frac{1}{\gamma_i}.$$
Notice that $a_{n+1}=a_n /\gamma_n.$
Thus, $a_n$ is increasing in $n$ if $0 < p \leq 1/2 $, and  decreasing in $n$ if $ 1/2<  p \leq 1. $ Hence, we have
\begin{eqnarray}\label{ewer}
\max_{ 2\leq i  \leq  n  }a_i = \left\{ \begin{array}{ll}
a_n \ \ \ & \textrm{if $0 <p \leq 1/2 $}\\
\\
(2p)^{-1} \ \ \ & \textrm{if $ 1/2 < p \leq 1.$}
\end{array} \right.
\end{eqnarray}
Define the filtration $\mathcal{F}_0=\{ \emptyset, \Omega\}$, $\mathcal{F}_n=\sigma\{\alpha_{i }, X_{1 }, \beta_{i },  Z_{i}  : 1\leq i\leq n\}, \, n\geq  1, $ and denote
\begin{eqnarray}
  a_nS_{n } - (2q-1)
 =   a_n \sum_{i=1}^n\alpha_{i }X_{\beta_{i }}(Z_{i}-1) +  a_nT_n  -(2q-1). \qquad  \label{mndef}
\end{eqnarray}
Let $$M_0=0,\ \ \ M_n=a_nT_n-(2q-1), \  \ n \geq 1.$$
As noticed in \cite{B18}, $\big(M_n, \mathcal{F}_n \big)_{n\geq0}$ is a martingale. Indeed, we have
$$\mathbf{E}[M_1| \mathcal{F}_0 ] =\mathbf{E} [a_{ 1}T_{ 1}-(2q-1) ] =0,$$ and for   $n\geq 1,$
\begin{eqnarray}
\mathbf{E}[M_{n+1} | \mathcal{F}_n] &=&\mathbf{E}[a_{n+1}T_{n+1}-(2q-1)| \mathcal{F}_n] \nonumber \\
&=&a_{n+1}T_n+a_{n+1}\mathbf{E}[\alpha_{n+1}] \mathbf{E}[X_{\beta_{n+1}}|\mathcal{F}_n]   -(2q-1) \nonumber\\
&=&a_{n+1} \bigg(T_n  + (2p-1)   \frac{T_n}{n} \bigg) -(2q-1) \nonumber\\
&=&a_nT_n -(2q-1)=M_n   \qquad a.s.  \nonumber
\end{eqnarray}
Moreover, we can rewrite $(M_n , \mathcal{F}_n)_{n\geq1}$ in the following additive form
\begin{equation}\label{fsdfsdf}
 M_n =\sum_{i=1}^{n} a_i\varepsilon_i :=  \sum_{i=1}^{n} \Delta M_{i} ,
\end{equation}
where $\varepsilon_i=T_i-\gamma_{i-1}T_{i-1}$ with the convention    $\gamma_0T_0=2q-1.$
Note that  $(\Delta M_{k}, \mathcal{F}_k)_{  k \geq 0 }$ is a sequence of martingale differences such that
\[
{\rm Var} (M_n) = \sum_{i=1}^n a_i^2 - (2p-1)^2 {\mathbf E}   \zeta_n,\ \  \text{ where }\ \zeta_n = \sum_{k=1}^{n-1} a^2_{k+1} \Big (  \frac{ T_k}{k} \Big )^2
\]
(see (A.7) in \cite{B18}).  In particular
\begin{equation} \label{boundmoment2}
{\rm Var} (M_n) \leq v_n \ \ \   \text{ and }\ \ \  \Vert a_n T_n \Vert_2^2 \leq 2 ( 1 + v_n) .
\end{equation}
The following lemma shows that $\Delta M_{ k}$ is bounded.
\begin{lemma}\label{lem1}
	For  $k\geq1$ and $p \in (0, 1]$, it holds $|\Delta M_k|\leq 2a_k$.
\end{lemma}
\noindent\emph{Proof.}
We have  $$|\Delta M_1|=|a_1T_1-(2q-1)| \leq  |X_1|+ |2q-1| \leq 2=2a_1,$$ and, for $k\geq 2,$
\begin{eqnarray*}
	\Delta M_k = a_k T_k-a_{k-1}T_{k-1}
	 = a_kX_k-a_k\frac{T_{k-1}}{k-1} (2p-1).
\end{eqnarray*}
Since $|X_k| \leq 1$, we have   $|T_{k }|\leq k $.  Hence, for $k\geq2,$
$$|\Delta M_k|\leq a_k|X_k|+a_k|2p-1| \leq 2a_k.$$
This completes the proof of the lemma.\hfill\qed

\smallskip

Next lemma gives estimates of the conditional moments of order $3$ and $4$ of the sequence of martingale differences $(\Delta M_{ k}, \mathcal{F}_{k})_{  k \geq 1 }$.
 \begin{lemma}\label{lem2}
 For  $k\geq 1$, the following equalities hold:
\[\mathbf{E}[  \Delta M_{ k}^3 | \mathcal{F}_{k-1}] =   2 (2p-1)a^3_k\bigg( -\frac{T_{k-1}}{k -1 }  + (2p-1)^2 \frac{T_{k-1}^3}{(k -1)^3 }   \bigg) \]
 and
\[ \mathbf{E}[  \Delta M_{ k}^4 | \mathcal{F}_{k-1}]    =  a_k^4  +   a_k^4 (2p-1)^2 \bigg( 2 \frac{T^2_{k-1}}{(k -1)^2 }  -3  (2p-1)^2 \frac{T_{k-1}^4}{(k -1)^4 }   \bigg),
 \]
 with  the convention  $ T_0/0=1.$ In particular,
 \[ \Big|\mathbf{E}[  \Delta M_{ k}^3 | \mathcal{F}_{k-1}] \Big| \ \leq \  4 a^3_k\bigg|\frac{T_{k-1}}{k -1 }   \bigg|   \quad \text{and}  \quad \Big|\mathbf{E}[  \Delta M_{ k}^4 | \mathcal{F}_{k-1}]  - a_k^4 \Big| \ \leq \  5 a^4_k \frac{T^2_{k-1}}{(k -1)^2 } .   \]
\end{lemma}
\noindent\emph{Proof.}  The equality for the conditional moment of order $4$ comes from \cite[(A.5)]{B18}. It remains to prove the equality for
 the conditional moment of order $3$.  Note first that
\begin{eqnarray}
 &&\mathbf{E}[  \Delta M_{ k}^3 | \mathcal{F}_{k-1}] = a^3_k  \mathbf{E}[ (T_k- \gamma_{k-1}T_{k-1})^3 | \mathcal{F}_{k-1}] \nonumber \\
 &&\ \ \ \ \ = \ a^3_k \Big( \mathbf{E}[  T_k ^3 | \mathcal{F}_{k-1}] - 3\gamma_{k-1}T_{k-1}\mathbf{E}[  T_k ^2 | \mathcal{F}_{k-1}]  +3\gamma_{k-1}^2T_{k-1}^2\mathbf{E}[  T_k   | \mathcal{F}_{k-1}] - \gamma_{k-1}^3 T_{k-1} ^3    \Big). \nonumber
\end{eqnarray}
For $\mathbf{E}[  T_k ^3 | \mathcal{F}_{k-1}]$, we have the following estimation
\begin{eqnarray}
 \mathbf{E}[  T_k ^3 | \mathcal{F}_{k-1}]&=& \mathbf{E}[  (T_{k-1}+ \alpha_{k}X_{\beta_k} )^3 | \mathcal{F}_{k-1}] \nonumber \\
 &=&   T_{k-1} ^3 + 3 T_{k-1} ^2  \mathbf{E}[ \alpha_{k}X_{\beta_k} | \mathcal{F}_{k-1}]
  \nonumber \\
  && \quad \quad \quad
  + 3 T_{k-1}    \mathbf{E}[ (\alpha_{k}X_{\beta_k})^2 | \mathcal{F}_{k-1}] + \mathbf{E}[   ( \alpha_{k}X_{\beta_k} ) ^3 | \mathcal{F}_{k-1}]\nonumber \\
 &=&    T_{k-1} ^3  + 3(2p-1) \frac{T_{k-1} ^3}{k-1}   + 3 T_{k-1}   + (2p-1) \frac{T_{k-1}}{k-1}  \nonumber \\
  &=&  (3 \gamma_{k-1} -2)  T_{k-1} ^3    +  (  \gamma_{k-1} +2)  T_{k-1}   .  \nonumber
\end{eqnarray}
Notice that $\mathbf{E}[T_k^2|  \mathcal{F}_{k-1}]  =    (2 \gamma_{k-1} -1) T_{k-1}^2    + 1$ and $\mathbf{E}[ T_k |  \mathcal{F}_{k-1}]  =  \gamma_{k-1}   T_{k-1}$ (see the relations (A.3) and (2.3) in \cite{B18}).
Hence, we get
\begin{eqnarray}
 \mathbf{E}[  \Delta M_{ k}^3 | \mathcal{F}_{k-1}]  =  2a^3_k \Big(  T_{k-1}(1 -\gamma_{k-1}) +  T_{k-1}^3 (\gamma_{k-1} -1)^3     \Big). \nonumber
\end{eqnarray}
The first equality follows from the fact that  $\gamma _n=1+\frac{2p-1}{n}$.
To get the  inequalities, it suffices to notice that  $|T_{k }|\leq k $. \hfill\qed

\section{Proof of Theorem  \ref{th3.1}}
\setcounter{equation}{0}
We first give a proof of item [i].
For any $n \geq 1$, let
\begin{eqnarray} \label{kndssl}
H_n=\sum_{i=1}^n\alpha_{i }X_{\beta_{i}}(Z_{i}-1) .
\end{eqnarray}
By  \eqref{mndef}, recall that  $
S_{n }  =  H_n +T_n$.
Let ${\mathcal G}_n= \sigma(\alpha_i, \beta_i, X_1, Z_i, 1\leq i\leq n).$  Note that $(H_n,{\mathcal G}_n)_{n \geq1}$ is a zero-mean, square integrable martingale such that
 $\textrm{Var}(\alpha_{i }X_{\beta_{i}}(Z_{i}-1))=\sigma^2$ and $ |\alpha_{i }X_{\beta_{i}}(Z_{i}-1)| =  |Z_{i}-1|$. Then, by \cite[Corollary 4.2]{HH80},   \begin{eqnarray}\label{knsl}
  \limsup_{n\rightarrow \infty}\frac{  |H_n| }{\sqrt{n \log \log n}  }  = \sqrt{2} \sigma \ \ \ a.s.
\end{eqnarray}
On another hand, by Theorem 3.2 of Bercu \cite{B18},
$$ \limsup_{n\rightarrow \infty}\frac{  |T_n| }{\sqrt{n \log \log n}  }  =\sqrt{\frac{2}{3-4p}} \ \ \ a.s.  $$
Therefore
\begin{eqnarray*}
\limsup_{n \rightarrow \infty} \frac{|S_{n}|}{\sqrt{n \log \log n} } &\leq& \limsup_{n \rightarrow \infty} \frac{|H_{n}|}{\sqrt{n\log \log n} } + \limsup_{n \rightarrow \infty} \frac{|T_{n}|}{ \sqrt{n\log \log n} }\\
 &\leq& \sqrt{2} \sigma + \sqrt{\frac{2}{3-4p}}\ \ \ \ a.s., \
\end{eqnarray*}
which gives \eref{fsv12}.

We turn to the proof of item [ii]. By Theorem 3.5 of Bercu \cite{B18}, we have
$$ \limsup_{n\rightarrow \infty}\frac{  |T_n| }{\sqrt{n  \log n \log \log \log n }  }  =\sqrt{2 } \ \ \ a.s.  $$
Since $S_n = H_n +T_n$, we get  from \eref{knsl} that
$$\limsup_{n \rightarrow \infty} \Bigg| \frac{|S_{n}|- |T_n|}{\sqrt{n \log n \log \log \log n } } \Bigg|\leq\limsup_{n \rightarrow \infty} \frac{|H_{n}|}{\sqrt{n  \log n \log \log \log n }  }     =  0\ \ \ \ a.s.,  $$
and  \eqref{fsv13} follows.

We complete the proof by proving item  [iii] that is when  $ 3/4 <p \leq 1$.  Since in this case $ n^{2p-1}  / \sqrt{n} \rightarrow 0$,   \eref{knsl} entails that
$ \lim_{n \rightarrow \infty} \frac{H_{n}}{n^{ 2p-1} } =0$. Hence the result follows from Theorems 3.7 and 3.8 in Bercu \cite{B18} and the fact that $S_n =H_n + T_n$.


\section{Proof of Theorem  \ref{th0}}  \label{secp}
\setcounter{equation}{0}
By \eref{kndssl}, we have
\begin{equation}\label{dec1}
\frac{a_nS_n -(2q-1)}{ \sqrt{v_n+ n a_n^2 \sigma^2}} =\frac{a_nH_n}{\sqrt{v_n+ n a_n^2 \sigma^2}} +\frac{M_n}{\sqrt{v_n+ n a_n^2 \sigma^2}}=:U_n+V_n ,
\end{equation}
where $H_n$ is defined in \eqref{kndssl}.
Denote $\mathcal{F}=\sigma(\alpha_i, \beta_i, X_1, i\geq 1).$ For $s, t \in \mathbf{R},$ write
\[f_{n}(s,t) =\mathbf{E} \exp\{(itU_n+isV_n)\} =\mathbf{E}\big[\mathbf{E}[\exp\{itU_n\}|\mathcal{F}\,] \exp\{isV_n\} \big] , \]
the joint characteristic function of $(U_n, V_n)$.   Since $(Z_i)_{i \geq 1}$ is a sequence of i.i.d. r.v.'s,
\[\lim_{n\rightarrow \infty} \mathbf{E}[\exp\{itU_n\}|\mathcal{F}\,] =\varphi(\sigma_1t)\ \ \ \ \textrm{a.s.,}\ \ \
\textrm{where} \ \ \    \sigma^2_1= \lim_{n\rightarrow \infty} \frac{na_n^2 \sigma^2}{v_n + na_n^2\sigma^2 } \] and $\varphi(t) = \exp\{- t^2/2\}$.
 Now
 $$f_n(s,t) =\mathbf{E}\Big[\big(\textbf{E}[\exp\{itU_n\}|\mathcal{F}\,]-\varphi(\sigma_1t) \big) \exp\{isV_n\} \Big] +\varphi(\sigma_1 t)\mathbf{E}[\exp\{isV_n\}].$$
From the CLT for the usual ERW (see Theorems 3.3 and 3.4 in Bercu \cite{B18}), we know that $$\lim_{n\rightarrow \infty} \mathbf{E}[\exp\{isV_n\}] =\varphi(\sigma_2s),\ \ \  \textrm{where} \ \ \ \sigma_2^2= \lim_{n\rightarrow \infty} \frac{ v_n}{v_n +na^2_n \sigma^2 }.$$
 Consequently, it holds $$\lim_{n\rightarrow \infty}f_n(s,t) =\varphi(\sigma_1t)\varphi( \sigma_2s).$$
In particular, since $\sigma^2_1+\sigma_2^2= 1$, we get $\lim_{n\rightarrow \infty}f_n(t,t) =\varphi( t)$, which implies that $(U_n+V_n)_{n\geq1}$ converges in distribution to a standard normal random variable.

\section{Proof of Theorem  \ref{thmBEERW}}
\setcounter{equation}{0}

The proof of  Theorem  \ref{thmBEERW} is a refinement on the argument of Lemma 3.3 of Grama and Haeusler \cite{GH00},
 where the authors obtained the best possible Berry-Esseen bound for martingales with bounded differences.
Compared to  \cite{GH00},   the main challenge of our proof comes from
the fact that the conditional variance of martingales does not converge to a constant in $L^\infty$-norm.
Our proof will be based on Lindeberg's telescoping sums argument, see Bolthausen \cite{B82}. The following two  technical lemmas  due to Bolthausen   \cite{B82} will be needed.
\begin{lemma}\cite[Lemma 1]{B82}
\label{LEMMA-APX-1}Let $X$ and $Y$ be random variables. Then
\[
\sup_{u \in \mathbf{R}}\Big| \mathbf{P}\left( X\leq u\right) -\Phi \left( u \right) \Big|
\leq c_1\sup_{u \in \mathbf{R}}\Big| \mathbf{P}\left( X+Y\leq u\right) -\Phi \left( u \right)
\Big| +c_2 \big\| \mathbf{E}\left[ Y^2|X\right] \big\| _\infty ^{1/2},
\]
where $c_1=2$ and $c_2 = 5 / \sqrt{2 \pi}$.
\end{lemma}
\begin{lemma} \cite[Lemma 2]{B82}
\label{LEMMA-APX-2}Let $G$ be an integrable function on $\mathbf{\mathbf{R}}$ of bounded variation $\|G\|_V$,
$X$ be a random variable and $a,$ $b\neq 0$ be real numbers. Then
\[
\mathbf{E}\bigg[ \, G\left( \frac{X+a}b\right) \bigg] \leq \|G\|_V \sup_u\Big| \mathbf{P}\left( X\leq u\right)
-\Phi \left( u\right) \Big| + \|G\|_1 \, |b|,
\]
where  $\|G\|_1$ is the $L_1(\mathbf{R})$ norm of $G.$
\end{lemma}

We shall divide the proof in two steps according to the values of $p$.

\smallskip

\noindent \textit{1) Case $p \in (0,3/4)$.} Let us introduce some notations that will be used all along the proof. Let $\beta$ be a constant greater than one that will be specified later. Denote
\begin{equation} \label{defkappan}
\kappa_n^2=\max_{[n/2]\leq i \leq n} a_i^2 \, , \,   \delta_n^2 = \beta \kappa_n^2 / v_n  ,\end{equation}
and for $ 0 \leq k \leq n$,
\begin{equation} \label{deftkn}
t_{k,n}= \big ( \beta \kappa_n^2  +  v_n-v_k \big )^{1/2} \, , \, A_{k}  =  t^2_{k,n}  /v_n \, , \,  u_{k,n}= \beta  +  \beta^{-1} ( 1 + \beta) \sum_{i=k+1}^n  a_i^2 / \kappa_n^2    .
\end{equation}
Note that $A_k$ is  non-increasing in $k,$ and satisfies
$A_0=\delta_n^2+  1 $ and $\displaystyle A_n=\delta_n ^2  $.
Moreover, for  $u, x\in \mathbf{R}$ and $y > 0,$ set, for brevity,
\begin{equation}
\Phi _u(x,y)=\Phi \bigg( \frac{u-x}{\sqrt{y}}  \bigg) .  \label{RA-7}
\end{equation}
Let $\mathcal{N}_1$  be  a  standard normal random variable, which is also independent of $T_{n}$.
For $ 1 \leq k  \leq n,$ denote
$$ \ \hat{K}_0=0, \ \   \ \  \hat{K}_k =\frac{1}{\sqrt{v_n}} \Big(a_kT_k  -(2q-1) \Big),\ \ \ \ \xi_k=\hat{K}_k- \hat{K}_{k-1}   \ \ \ \    \textrm{and } \ \ \ \ \   \sigma_k^2= \frac{a_k^2}{v_n } .$$
Then we have
 $$\frac{a_nT_n -(2q-1)}{ \sqrt{v_n }}\ = \  \hat{K}_n  = \sum_{k=1}^n \xi_k .$$
Next, we estimate the upper bound of $\mathbf{K}(\hat{K}_n ).$ By Lemma \ref{LEMMA-APX-1}, we get
\[
 \mathbf{K}(\hat{K}_n )
\leq c_1\sup_{u \in \mathbf{R}}\bigg|
 \mathbf{P} \Big(   \hat{K}_n +  \delta_n\mathcal{N}_1  \leq u \Big)  -  \Phi (u)\bigg|
+c_2 \delta_n  \nonumber  \\
=  c_1\sup_{u \in \mathbf{R}}\bigg|
\mathbf{E} [\Phi _u(\hat{K}_n,A_n)]- \Phi (u)\bigg|
+c_2\, \delta_n .
\]
Hence
\[
  \mathbf{K}(\hat{K}_n )
\leq  c_1\sup_{u \in \mathbf{R}}\bigg| \mathbf{E} [\Phi _u(\hat{K}_n,A_n)]-\mathbf{E} [\Phi _u(\hat{K}_0,A_0)]\bigg|  \  +\, c_1\sup_{u \in \mathbf{R}}\bigg| \Phi \Big(\frac{u}{\sqrt{\delta_n ^2+ 1}}  \Big)-\Phi (u)\bigg| +c_2\, \delta_n  . \]
There exist positive constants $c_3$ and $c_4$ such that for any $n \geq 1$,
\begin{eqnarray}
 \bigg| \Phi \Big(\frac{u}{\sqrt{\delta_n ^2+ 1}}  \Big)-\Phi (u)\bigg| \leq  c_3 \bigg|\frac{1}{\sqrt{\delta_n ^2+ 1}}  -1 \bigg|
\leq c_4\delta_n ^2 .
\end{eqnarray}
Therefore, since $\delta_n^2 \leq \beta$,
\begin{eqnarray}
\mathbf{K}(\hat{K}_n ) \ \leq \  c_1\sup_{u \in \mathbf{R}}\Big| \mathbf{E} [\Phi _u(\hat{K}_n,A_n)]-\mathbf{E} [\Phi _u(\hat{K}_0,A_0)]\Big|   +c_5\delta_n .   \label{dsdff}
\end{eqnarray}
Next, we give an estimation for $\mathbf{E} [\Phi _u(\hat{K}_n,A_n)]-\mathbf{E} [\Phi _u(\hat{K}_0,A_0)]$.
Write first
\[
\mathbf{E} [\Phi _u(\hat{K}_n,A_n)]-\mathbf{E} [\Phi _u(\hat{K}_0,A_0)]= \sum_{k=1}^n\mathbf{E}
\big[ \Phi _u(\hat{K}_k,A_k)-\Phi _u(\hat{K}_{k-1},A_{k-1})  \big] .
\]
Using the fact that
\[
\frac{\partial ^2}{\partial x^2}\Phi _u(x,y)=2\frac \partial {\partial
y}\Phi _u(x,y),
\]
and that $\mathbf{E} [ \xi_k | {\mathcal F}_{k-1} ] =0$,
we obtain
\begin{equation}\nonumber
\mathbf{E}[ \Phi _u(\hat{K}_n,A_n)]-\mathbf{E} [\Phi _u(\hat{K}_0,A_0)]=I_1+I_2-I_3,
\end{equation}
where
\begin{eqnarray}
I_1&=& \sum_{k=1}^n \mathbf{E} \bigg[  \frac{}{} \Phi _u(\hat{K}_k,A_k)-\Phi
_u(\hat{K}_{k-1},A_k)  \nonumber \\
&& \ \ \ \ \ \ \ \ \ \ \    -\, \frac \partial {\partial x}\Phi
_u(\hat{K}_{k-1},A_k)\xi_k-\frac 12\frac{\partial ^2}{\partial x^2}\Phi
_u(\hat{K}_{k-1},A_k)\xi _k^2   \bigg] ,   \nonumber\\
I_2 &=& \frac 12 \sum_{k=1}^n \mathbf{E}   \bigg[\frac{\partial ^2}{\partial x^2}\Phi
_u(\hat{K}_{k-1},A_k)\Big( \mathbf{E} [ \xi _k^2   | \mathcal{F}_{k-1}]   - \sigma_k^2\Big) \bigg] ,\quad \quad \
 \nonumber \\
I_3&=&\sum_{k=1}^n\mathbf{E}   \bigg[  \Phi _u(\hat{K}_{k-1},A_{k-1})-\Phi
_u(\hat{K}_{k-1},A_k)-\frac \partial {\partial y}\Phi _u(\hat{K}_{k-1},A_k)
\sigma_k^2 \bigg].  \nonumber
\end{eqnarray}
From (\ref{dsdff}), we deduce that
\begin{eqnarray}\label{inea77}
  \mathbf{K}(\hat{K}_n )\  \leq \ C\, \Big(  |I_1|+|I_2|+|I_3| + \delta_n \Big).
\end{eqnarray}
In the sequel, we  give some estimates for $I_1,$ $I_2$ and $I_3.$ The notation
$\vartheta_i$ stands for some values or random variables satisfying $0\leq \vartheta_i  \leq 1$.

\emph{\textbf{a)} Control of} $I_1.$
To shorten notations, denote
$\displaystyle  \hat{H}_{k-1} (u)  = \frac{u-\hat{K}_{k-1}}{\sqrt{A_k}} .$ For $1\leq k \leq n,$ we have
\begin{align*}
R_k & := \Phi _u(\hat{K}_{k},A_k)-\Phi
_u(\hat{K}_{k-1},A_k)    -\, \frac \partial {\partial x}\Phi
_u(\hat{K}_{k-1},A_k)\xi_k-\frac 12\frac{\partial ^2}{\partial x^2}\Phi
_u(\hat{K}_{k-1},A_k)\xi _k^2  \\
 & =    \Phi \Big( \hat{H}_{k-1} (u)  - \frac{\xi_k}{\sqrt{A_k}}\Big) -\Phi( \hat{H}_{k-1} (u) )
    +\, \Phi'(\hat{H}_{k-1} (u) )\frac{\xi_k}{\sqrt{A_k}}   \\ & \quad \quad \quad \quad  \quad \quad   -\frac 12\Phi''(\hat{H}_{k-1} (u) )\Big(\frac{\xi_k}{\sqrt{A_k}}\Big)^2 .
\end{align*}
By  Taylor expansion at order $4$,   we deduce that
\begin{eqnarray}\label{gd55}
I_1  & =& \sum_{k=1}^n \mathbf{E}  R_{k,1} +\sum_{k=1}^n   \mathbf{E} R_{k,2},
\end{eqnarray}
where
\begin{eqnarray}\nonumber
 R_{k,1}=  \,-\frac16 \Phi'''\big(\hat{H}_{k-1} (u)   \big)  \Big(\frac{\xi_k}{\sqrt{A_k}}\Big)^3  \ \  \ \textrm{and} \ \ \  R_{k,2}= \,\frac1{24} \Phi^{(4)}\Big(\hat{H}_{k-1} (u)  - \vartheta_k  \frac{\xi_k}{\sqrt{A_k}}\Big)  \Big(\frac{\xi_k}{\sqrt{A_k}}\Big)^4 .
\end{eqnarray}
 Next, we handle the term $| \mathbf{E}  R_{k,1} |.$ As $\hat{H}_{k-1} (u) $ is $\mathcal{F}_{k-1}$-measurable, we have
\[
 \big| \mathbf{E}  R_{k,1}\big| \leq    \frac{1}{A_k^{3/2}} \mathbf{E} \Big[ \Big|\Phi'''\big(\hat{H}_{k-1} (u)   \big) \Big| \big| \mathbf{E}[ \xi_k^3| \mathcal{F}_{k-1} ] \big| \Big] \, \leq \, C \frac{1}{A_k^{3/2}} \mathbf{E}    \big| \mathbf{E}[ \xi_k^3| \mathcal{F}_{k-1} ] \big|   .
\]
 By Lemma \ref{lem2}, we have $ \big|\mathbf{E}   \xi_{1}^{3}   \big|\leq4 a_1^3/v_n^{3/2},$ and for  $2 \leq k \leq n,$
\[
 \mathbf{E}   \big|\mathbf{E}  \big[ \xi_{k}^{3} \big| \mathcal{F}_{k-1} \big] \big| \ll     \Big( \frac{a_k}{ \sqrt{v_n}}\Big)^{3}  \mathbf{E}    \bigg|\frac{ T_{k-1} }{k-1}   \bigg| \ll      \frac{a_k^3 (\mathbf{E}    \big| a_{k-1} T_{k-1}    \big|^2)^{1/2}}{a_{k-1} v_n^{3/2} (k-1)}  . \]
Taking into account \eqref{boundmoment2}, it follows that for  $2 \leq k \leq n,$
\[
 \mathbf{E}   \big|\mathbf{E}  \big[ \xi_{k}^{3} \big| \mathcal{F}_{k-1} \big] \big| \ll      \frac{a_k^3 ( 1+ \sqrt{v_{k-1}} )}{a_{k-1} v_n^{3/2} (k-1)} .
\]
Hence,  we get
\begin{equation} \label{termeRk1-1}
\sum_{k=1}^n  \big| \mathbf{E}  R_{k,1}\big|   \ll  \frac{1}{v_n^{3/2}}+  \sum_{k=2}^n\frac{1}{A_k^{3/2}}   \frac{a_k^3 \sqrt{v_{k-1}} }{a_{k-1} v_n^{3/2} (k-1)} .
\end{equation}
Next we give an estimation for the last bound. By (\ref{a10}) and (\ref{a15}), and since $p \in (0, 3/4)$, $ \displaystyle  \frac{  \sqrt{v_{k-1}}}{ a_{k-1} (k-1)  } =O\Big(\frac{ 1 }{\sqrt{k}} \Big) , k \rightarrow \infty$. Note also that $A_k v_n= t_{k,n}^2$ and, since $p \in (0,3/4)$,
$\max_{1 \leq k \leq [n/2]}  t_{k,n}^{-2} =  t_{[n/2],n}^{-2}  \ll v^{-1}_n $. Therefore, we get
\begin{equation} \label{termeRk1-2}
\sum_{k=2}^n\frac{1}{A_k^{3/2}}   \frac{a_k^3 \sqrt{v_{k-1}} }{a_{k-1} v_n^{3/2} (k-1)}
 \ll   v_n^{-3/2}\!\! \sum_{k \in [1, n/2]}  \frac{a_k^3 }{\sqrt{k} }
   +  \frac{\kappa_n}{\sqrt{n} } \sum_{k \in(  n/2, n]}  t_{k,n}^{-3} a_k^2  .
\end{equation}
For the first term  on the right-hand side of  (\ref{termeRk1-2}), we infer that
\begin{equation} \label{termeRk1-3}
v_n^{-3/2}\!\! \sum_{k \in [1, n/2]}  \frac{a_k^3 }{\sqrt{k} } \ll v_n^{-3/2} {\bf 1}_{p >7/12} +  v_n^{-3/2}(\log n) {\bf 1}_{p =7/12}  + n^{-1} {\bf 1}_{p <7/12}  .
\end{equation}
To deal with the second term on the right-hand side of (\ref{termeRk1-2}), we  first notice that for any integer $k$ such that $ [n/2]  \leq k \leq n$,
\begin{equation} \label{ine52}
 \frac{a_k^2}{ t_{k,n}^2}  \leq \log \big (  u_{k-1,n}/u_{k,n}   \big )  \text{ and } \frac{\kappa_n^2 u_{k-1,n}}{3} \leq  t_{k,n}^2 \leq \kappa_n^2 u_{k,n}
\end{equation}
(see inequalities (5.2)  and (5.4) in \cite{DMR22}). Therefore, for any  $\alpha \geq 0 $,
\begin{multline*}
{\mathcal T}_n (\alpha):= \sum_{k=[n/2]+1}^n  t_{k,n}^{\alpha} \frac{a_k^2}{ t_{k,n}^2} \leq   \kappa_n^\alpha  \sum_{k=[n/2]+1}^n   u_{k,n}^{\alpha/2}  \int_{u_{k,n}}^{ u_{k-1,n}} \frac{1}{x} dx  \\
\leq   \kappa_n^\alpha  \sum_{k=[n/2]+1}^n      \int_{u_{k,n}}^{ u_{k-1,n}}  x^{\alpha/2-1}dx
\end{multline*}
and, for any $\alpha <0$,
\[
{\mathcal T}_n (\alpha) \leq   \kappa_n^\alpha 3^{-\alpha/2} \sum_{k=[n/2]+1}^n   u_{k-1,n}^{\alpha/2}  \int_{u_{k,n}}^{ u_{k-1,n}} \frac{1}{x} dx
\leq   \kappa_n^\alpha 3^{-\alpha/2} \sum_{k=[n/2]+1}^n      \int_{u_{k,n}}^{ u_{k-1,n}}  x^{\alpha/2-1}dx .
\]
So, overall, for any real $\alpha$, using the fact that $p \in (0,3/4)$, we get
\begin{align} \label{ine52cons}
{\mathcal T}_n (\alpha)= \sum_{k=[n/2]+1}^n  t_{k,n}^{\alpha} \frac{a_k^2}{ t_{k,n}^2} &  \leq     \kappa_n^\alpha  ( 1 + 3^{-\alpha/2} )    \int_{\beta }^{ u_{[n/2],n}}  x^{\alpha/2-1}dx
\nonumber \\ & \ll   \kappa_n^\alpha  \beta^{\alpha/2} {\mathbf 1}_{\alpha < 0}  +  ( \log n ) {\mathbf 1}_{\alpha = 0}   +  v_n^{\alpha/2} {\mathbf 1}_{\alpha > 0}  .
\end{align}
Therefore, by \eqref{ine52cons} applied with $\alpha =-1$, we get
\begin{equation} \label{termeRk1-4}
\frac{\kappa_n}{\sqrt{n} } \sum_{k \in(  n/2, n]}  t_{k,n}^{-3} a_k^2  \ll  n^{-1/2} .
\end{equation}
So, overall,  starting from \eqref{termeRk1-1} and taking into account \eqref{termeRk1-2}, \eqref{termeRk1-3}, \eqref{termeRk1-4} and that $(\log n) v_n^{-3/2}  {\bf 1}_{p=7/12}  \ll n^{-1 }\log n $, we derive that, for $p\in (0,3/4)$,
\begin{eqnarray}\label{fs17}
\sum_{k=1}^n  \big| \mathbf{E}  R_{k,1}\big|   \ll  \frac{1}{v_n^{3/2}}+  \frac{1}{ \sqrt{n}} .
\end{eqnarray}
Next, we give an estimate of   $ |\mathbf{E}  R_{k,2}|$. By Lemma \ref{lem1}, note that $\displaystyle
\Big|\frac{\xi_{k}}{ \sqrt{A_{k}}}\Big| \leq \frac{2a_k}{t_{k,n}} $. Hence, if $k \geq [n/2]+1$, $\displaystyle \Big|\frac{\xi_{k}}{ \sqrt{A_{k}}}\Big| \leq 2/\sqrt{\beta}$. Next, when $k \leq [n/2]$ and $p \in [0,1/2]$, we still have $\displaystyle \Big|\frac{\xi_{k}}{ \sqrt{A_{k}}}\Big| \leq 2/\sqrt{\beta}$. On another hand, when $k \leq [n/2]$ and $p >1/2$, we have
$\displaystyle \Big|\frac{\xi_{k}}{ \sqrt{A_{k}}}\Big| \leq  \min \big (  2/  \sqrt{ \beta \kappa_n^2} , 2/\sqrt{v_n - v_{[n/2]}} \big ) $. So, in each case, we can select $\beta $ large enough in such a way that
$\displaystyle \Big|\frac{\xi_{k}}{ \sqrt{A_{k}}}\Big|  \leq 1$. From now on, $\beta$ is selected this way. Thus we have
\begin{eqnarray*}
   |R_{k,2}|     & \leq& \frac1{24} \Big|\Phi^{(4)}\Big(\hat{H}_{k-1} (u)  - \vartheta_k  \frac{\xi_k}{\sqrt{A_k}}\Big) \Big| \bigg(\frac{\xi_k}{\sqrt{A_k}}\bigg)^4  \\
  &\leq& G(\hat{H}_{k-1} (u) ) \bigg(\frac{\xi_k}{\sqrt{A_k}}\bigg)^{4},
\end{eqnarray*}
where $G(z)=\sup_{|t-z|  \leq   1  } |\Phi^{(4)}\big(t  \big) |.$ By Lemma   \ref{lem1}, $ \xi_{k}^4 \leq 2^4 a_k^4 v_n^{-2}$. Therefore,
\begin{equation*}
\sum_{k=1}^n   |  \mathbf{E}   R_{k,2}  |  \ll   \sum_{k=1}^n \frac{\mathbf{E}[ G(\hat{H}_{k-1} (u) )]}{ A_{k}^{2} }   \bigg( \frac{a_k}{ \sqrt{v_n}}\bigg)^{4}   .
\end{equation*}
By the definition of
$ \displaystyle \hat{H}_{k-1} (u) $  and  Lemma \ref{LEMMA-APX-2}, it follows that there exists a positive constant $C$ such that for any $k \geq 1$ and any $u \in {\mathbf R}$,
\[
\mathbf{E}[ G(\hat{H}_{k-1} (u) )]  \leq   C\,  \mathbf{K}\big( \hat{K}_{k-1}    \big) + C \sqrt{A_k} .
\]
Now, by Lemma \ref{LEMMA-APX-1},
\[
\mathbf{K}\big( \hat{K}_{k-1}    \big)   \leq  C \mathbf{K}\big( \hat{K}_{n}    \big) + C \Vert  \mathbf{E}[   (  \hat{K}_{n}  -  \hat{K}_{k-1}  )^2  |  \hat{K}_{k-1}  ] \Vert_{\infty}^{1/2},
\]
and, by the martingale property of $(M_n, \mathcal{F}_{n})_{n\geq1}$ and Lemma \ref{lem1},
\[
 \mathbf{E}[   (  \hat{K}_{n}  -  \hat{K}_{k-1}  )^2  |  \hat{K}_{k-1}  ] = v_n^{-1}  \sum_{i=k}^n \mathbf{E}[   (  \Delta M_i )^2  |  \hat{K}_{k-1}  ]  \leq 4  v_n^{-1} ( v_n - v_{k-1}).
\]
Therefore, it holds 
\begin{equation} \label{conslmas7172}
\mathbf{E}[ G(\hat{H}_{k-1} (u) )]  \ll \mathbf{K}\big( \hat{K}_{n}    \big)  +  \sqrt{A_{k-1}}.
\end{equation}
So, overall, we have
\begin{equation} \label{B1forRk2}
\sum_{k=1}^n   \big |  \mathbf{E}  R_{k,2}      \big|
\ll \sum_{k=1}^n \frac{1}{ A_{k}^{2} }  \bigg( \frac{a_k}{ \sqrt{v_n}}\bigg)^{4} \mathbf{K}\big( \hat{K}_{n} \big) \  + \sum_{k=1}^n \frac{1}{ A_{k-1}^{3/2}  }  \bigg( \frac{a_k}{ \sqrt{v_n}}\bigg)^{4}  .
\end{equation}
Next, we give some estimations for the right-hand side of the last inequality. Recall that $\max_{1 \leq k \leq [n/2]}  (A_k v_n)^{-1} =  t_{[n/2],n}^{-2}  \ll v^{-1}_n $.   Therefore
\[ \sum_{k=1}^n \frac{1}{ A_{k}^{2} }  \bigg( \frac{a_k}{ \sqrt{v_n}}\bigg)^{4}  \ll \frac{1}{v_n^2} \sum_{k=1}^{[n/2]} a_k^4 + \kappa_n^2  \sum_{k=[n/2] +1}^n  t_{k,n}^{-2}\frac{a^2_k}{ t_{k,n}^2} .\]
Now $  \sum_{k=1}^n a_k^4 \ll {\mathbf 1}_{p>5/8 } + (\log n) {\mathbf 1}_{p=5/8 } + a_n^2 v_n {\mathbf 1}_{p< 5/8 } $.  On another hand, by taking into account
\eqref{ine52cons} with $\alpha =-2$,  we get
\[
 \sum_{k=[n/2] +1}^n  t_{k,n}^{-2}\frac{a^2_k}{ t_{k,n}^2} \ll     \kappa_n^ {-2}  \beta^{-1} .
\]
It follows that
\[
\sum_{k=1}^n \frac{1}{ A_{k}^{2} }  \bigg( \frac{a_k}{ \sqrt{v_n}}\bigg)^{4} \ll   \beta^{-1} +  n^{-1} + v_n^{-2} \log n.
\]
Using similar arguments, we infer that
\[
 \sum_{k=1}^n \frac{1}{ A_{k}^{3/2}  }  \bigg( \frac{a_k}{ \sqrt{v_n}}\bigg)^{4} \ll  v_n^{-2} \sum_{k=1}^n a_k^4 + v_n^{-1/2} \kappa_n \ll v_n^{-2} \big( {\mathbf 1}_{p>5/8 } + (\log n) {\mathbf 1}_{p=5/8 } \big) + n^{-1/2}.
\]
It follows that, for $n$ large enough,
\begin{equation} \label{dfsc}
 \sum_{k=1}^n   \big |  \mathbf{E}  R_{k,2}      \big|
\ll \beta^{-1}  \mathbf{K}\big( \hat{K}_{n} \big) \  + v_n^{-2}  {\mathbf 1}_{p>5/8 } + n^{-1/2} .
\end{equation}
Starting from  (\ref{gd55}) and taking into account \eqref{fs17} and \eqref{dfsc}, we derive that
\begin{equation}\label{controlI1}
|I_1| \ll  \beta^{-1}  \mathbf{K}\big( \hat{K}_{n} \big) \  + \frac{1}{v_n^{3/2}}+  \frac{1}{ \sqrt{n}} .
\end{equation}

\emph{\textbf{b)} Control of} $I_2.$ We have
\[
\left|
I_2\right| \leq \Bigg|\sum_{k=1}^n \frac 1{2 A_k}\mathbf{E}   \Big[  \varphi'(\hat{H}_{k-1} (u) ) \Big(\mathbf{E} [ \xi _k^2   | \mathcal{F}_{k-1}]   - \sigma^2_k \Big) \Big] \Bigg|\ll \sum_{k=1}^n \frac 1{ A_k} \mathbf{E}   \Big |   \mathbf{E} [ \xi _k^2   | \mathcal{F}_{k-1}]   - \sigma^2_k    \Big|.
\]
Clearly, we have  $ \mathbf{E} [ \xi _1^2   | \mathcal{F}_{0}]  = \sigma^2_1 $. On another hand, by relation (A.4) in \cite{B18}, we get, for $2 \leq k \leq n,$
\begin{equation} \label{ddfsds}
\mathbf{E}[ \Delta M_{ k} ^2   |\mathcal{F}_{k-1}]=a_k^2   - a_k^2( 2p -1)^2 \Big( \frac{T_{k-1}}{k-1} \Big)^2 ,
\end{equation}
so that
\begin{eqnarray*}
 \big |   \mathbf{E} [ \xi _k^2   | \mathcal{F}_{k-1}]   - \sigma^2_k    \big|
 \leq    ( 2p -1)^2  \frac{a^2_k}{v_n }  \Big ( \frac{T_{k-1}}{k-1} \Big )^2   .
\end{eqnarray*}
Thus, by \eqref{boundmoment2},
\begin{align*}
\left|
I_2\right|  &  \ll  \sum_{k=2}^n \frac 1{  A_k}  \frac{a^2_k}{v_n a_{k-1}^2(k-1)^2 }   \mathbf{E} (a_{k-1}T_{k-1})  ^2
 \\
& \ll  \sum_{k \in [2, n/2]}  \frac 1{  A_k}  \frac{a^2_k}{v_n}  \frac{  v_{k-1}}{ a_{k-1}^2(k-1)^2 } + \sum_{k \in (n/2, n]}  \frac 1{  A_k}  \frac{a^2_k}{v_n }  \frac{  v_{k-1}}{ a_{k-1}^2(k-1)^2 } .
\end{align*}
By  (\ref{a10}) and (\ref{a15}), we have for $p \in (0, 3/4),$ $ \displaystyle  \frac{  v_{k-1}}{ a_{k-1}^2(k-1)^2 } \asymp \frac{1  }{k}, k\rightarrow \infty$, and
$$\max_{1 \leq k \leq [n/2]} (A_k v_n)^{-1} \ll v^{-1}_n .$$ Hence
\[
\sum_{k \in [2, n/2]}  \frac 1{  A_k}  \frac{a^2_k}{v_n}  \frac{  v_{k-1}}{ a_{k-1}^2(k-1)^2 }  \ll v_n^{-1} \sum_{k \in [2, n/2]}   \frac{a_k^2}{ k}  \ll n^{-1} + v_n^{-1}
+ n^{-1}  ( \log n  ) {\mathbf 1}_{p=1/2} .
\]
On another hand, by  \eqref{ine52cons} with $\alpha =0$,  we get
\[
 \sum_{k \in (n/2, n]}  \frac 1{  A_k}  \frac{a^2_k}{v_n }  \frac{  v_{k-1}}{ a_{k-1}^2(k-1)^2 } \ll  n^{-1} \sum_{k \in (n/2, n]}  \frac {a_k^2}{  t_{k,n}^2}  \ll  n^{-1} \log n .
\]
So, overall,
\begin{equation}\label{controlI2}
|I_2| \ll   \frac{\log n}{n}  + \frac{1}{v_n} .
\end{equation}

\emph{\textbf{c)} Control of} $I_3.$ Note first that
\[
\Phi _u(\hat{K}_{k-1},A_{k-1})-\Phi
_u(\hat{K}_{k-1},A_k)-\frac \partial {\partial y}\Phi _u(\hat{K}_{k-1},A_k)
\sigma_k^2 = \frac{1}{2} \frac{ \partial^2}{\partial y^2}\Phi _u(\hat{K}_{k-1},A_k + \vartheta _k \sigma^2_k)
\sigma_k^4 ,
\]
for some $\vartheta_k \in [0,1]$. Now,
\[
 \frac{ \partial^2}{\partial y^2}\Phi _u(x,z) = \frac{3}{4} \frac{u-x}{z^{5/2}} \varphi \big ( \frac{u-x}{\sqrt{z}} \big )  + \frac{1}{4} \frac{(u-x)^2}{z^{3}} \varphi' \big ( \frac{u-x}{\sqrt{z}} \big )
\]
and,  for any $\vartheta \in [0,1]$, $A_{k-1}/5 \leq A_k +\vartheta  \sigma^2_k \leq A_{k-1} $. The right-hand side of this inequality is trivial. To prove the left hand side, we note that
 $a_k \leq 2  a_{k+1}$ which entails that $a_k^2 \leq 4 \sum_{\ell =k+1}^n a^2_\ell$. Therefore, $v_nA_{k-1}  \leq  \beta \kappa_n^2 + 5 \sum_{\ell =k+1}^n a^2_\ell \leq 5 v_nA_{k} $.
 Let
$
g(t) =  ( | t|  + |t|^3)  \varphi (t)
$. For $t = (u-\hat{K}_{k-1}) / \sqrt{A_k  + \vartheta  \sigma^2_k}$ and $ \widetilde{H}_{k-1} (u) =  (u-\hat{K}_{k-1}) / \sqrt{A_{k-1} }$, we then get
\[
g(t) \leq \Big ( \sqrt{5} |  \widetilde{H}_{k-1} (u)   | +   5^{3/2}|  \widetilde{H}_{k-1} (u)   |^3  \Big )  \varphi ( \widetilde{H}_{k-1} (u) ) := \widehat{G}( \widetilde{H}_{k-1} (u) ) .
\]
So, overall,
\[
\left| I_3\right| \leq     \sum_{k=1}^n\frac 1{A_k^2} \,\mathbf{E}[  \widehat{G} (\widetilde{H}_{k-1} (u)
  )]\, \sigma_k^4 .
\]
Proceeding as to get \eqref{conslmas7172}, we infer that $\mathbf{E}[  \widehat{G} (\widetilde{H}_{k-1} (u)
  )] \ll  \mathbf{K}\big( \hat{K}_{n}    \big)  +  \sqrt{A_{k-1}} $. Therefore
  \[
  \left| I_3\right|   \ll  \sum_{k=1}^n\frac 1{A_k^2} \sigma_k^4 \  \mathbf{K}\big(  \hat{K}_{n}   \big ) + \sum_{k=1}^n\frac 1{A_k^{ 3/2}} \sigma_k^4  .
  \]
  Hence the right-hand side of the above inequality is the same as the one of inequality \eqref{B1forRk2}. Therefore, according to \eqref{dfsc}, we get, for $n$ large enough, that \begin{equation} \label{controlI3}
|I_3|
\ll \beta^{-1}  \mathbf{K}\big( \hat{K}_{n} \big) \  + v_n^{-2}  {\mathbf 1}_{p>5/8 } + n^{-1/2} .
\end{equation}
Starting from \eqref{inea77} and taking into account the estimates \eqref{controlI1}, \eqref{controlI2}, \eqref{controlI3} and the fact that $\delta_n \ll v_n^{-1} + n^{-1} $, it follows that there exists a positive constant $C$ (not depending on $n$) such that, for $n$ large enough,
\[
 ( 1 -  C \beta^{-1}  )  \mathbf{K}\big( \hat{K}_{n} \big) \ll  n^{-1/2}  + v_n^{-1}  ,
\]
and  Theorem  \ref{thmBEERW} (for $p \in (0,3/4)$)   follows by taking $\beta$ large enough so that $1 -  C \beta^{-1}   \geq 1/2$.

\smallskip

\noindent \textit{2) Case $p =3/4$.}  Recall that in this case, $v_n /(\log n) \rightarrow \pi/4$. Compared to the previous case, the differences are as follows. First we fix a $\gamma \in (0,1)$ and we define
\begin{equation} \label{defkappanp=3/4}
\kappa_n^2=\max_{[n^{\gamma}]\leq i \leq n} a_i^2 = a_{[n^{\gamma}]}^2 .
\end{equation}
Next, $ \delta_n^2 $, $t_{k,n}$, $A_k$ and $u_{k,n}$ are still defined as in \eqref{defkappan} and \eqref{deftkn}. It follows that
\[
\max_{1 \leq k \leq [n^{\gamma}]} (A_kv_n)^{-1} = \max_{1 \leq k \leq [n^{\gamma}]} t^{-2}_{k,n} = t^{-2}_{[n^\gamma],n} \ll v^{-1}_n \, .
\]
On another hand,
for any real $\alpha$,
\begin{align} \label{ine52consp=3/4}
{\mathcal T}_n (\alpha)= \sum_{k=[n^{\gamma}]+1}^n  t_{k,n}^{\alpha} \frac{a_k^2}{ t_{k,n}^2} & \leq     \kappa_n^\alpha  ( 1 + 3^{-\alpha/2} )    \int_{\beta }^{ t_{[n^{\gamma}],n}}  x^{\alpha/2-1}dx \nonumber \\ & \ll   \kappa_n^\alpha  \beta^{\alpha/2} {\mathbf 1}_{\alpha < 0}  +  (\log \log n ) {\mathbf 1}_{\alpha = 0}   +  v_n^{\alpha/2} {\mathbf 1}_{\alpha > 0}  .
\end{align}
Moreover we shall use in case $p=3/4$ that $ \displaystyle  \frac{  \sqrt{v_{k-1}}}{ a_{k-1} (k-1)  } =O\Big(\frac{ \sqrt{\log k} }{\sqrt{k}} \Big)  , k \rightarrow \infty$.  Taking into account all the differences pointed above, separating the sums from $1$ to $n$ into a sum from $1$ to $[n^\gamma]$  plus a sum from  $[n^\gamma]+1$ to $n$, and proceeding as for the case $p \in (0,3/4)$ we infer that
\[
|I_1| + |I_2|  + |I_3|
\ll \beta^{-1}  \mathbf{K}\big( \hat{K}_{n} \big)  + v_n^{-1} ,
\]
which completes the proof of Theorem  \ref{thmBEERW} in case $p=3/4$.

\section{Proof of Theorem  \ref{thm6}}
\setcounter{equation}{0}

Now we are in the position to prove Theorem \ref{thm6}.
Set any integer $k$ in $ [1,n]$, set $\tau_k^2=   na_n^2 \sigma^2 +v_k$.
For $ 1 \leq k  \leq n,$ denote  \[ \widetilde{H}_k=\frac{ 1}{\tau_n }\sum_{i=1}^k\alpha_{i  }X_{\beta_{i  }}(Z_{i}-1) \mbox{ and } \widetilde{K}_k =\frac{1}{\tau_n} \Big(a_kT_k  -(2q-1) \Big) .\]
Then we have
 $$\frac{a_nS_n -(2q-1)}{ \sqrt{v_n+ n a_n^2 \sigma^2}}\ = \ a_n \widetilde{H}_n    + \widetilde{K}_n  .$$
By the triangle inequality, we get
\begin{eqnarray}
 \mathbf{K}\bigg(  \frac{a_nS_n -(2q-1)}{ \sqrt{v_n+ n a_n^2 \sigma^2}} \bigg)
  \leq  P_1  +   P_2, \label{fsdop}
\end{eqnarray}
where
$$P_1= \sup_{u \in \mathbf{R}}\bigg|
 \mathbf{P} \bigg( a_n \widetilde{H}_n    + \widetilde{K}_n   \leq u \bigg)  -  \mathbf{P} \bigg( \mathcal{N}_1 \frac{ \sqrt{n}a_n \sigma }{\tau_n }     + \widetilde{K}_n  \leq u \bigg)\bigg| $$
  and
  $$P_2= \sup_{u \in \mathbf{R}} \bigg| \mathbf{P} \bigg( \mathcal{N}_1 \frac{\sqrt{n}a_n \sigma }{\tau_n }    + \widetilde{K}_n  \leq u \bigg) -  \Phi (u) \bigg |,$$
where $ \mathcal{N}_1$ is a standard normal random variable independent of  $( \alpha_{k }, \beta_{k},  X_{k}, Z_k,   k \geq  0  )$. We first give an upper bound of $P_1.$
Let $\mathcal{F}= \sigma \{\alpha_{k }, \beta_{k},  X_{1},  k \geq  0  \}.$
By the classical Berry-Esseen bound for independent random variables (see \cite{E42}),
we have for   $ \rho \in (0, 1],$
\begin{eqnarray}
 && \bigg| \mathbf{P} \bigg( a_n \widetilde{H}_n    + \widetilde{K}_n   \leq u  \, \bigg| \,  \mathcal{F} \bigg)  -  \mathbf{P} \bigg( \mathcal{N}_1 \frac{\sqrt{n}a_n \sigma }{\tau_n }    + \widetilde{K}_n  \leq u \, \bigg| \, \mathcal{F}  \bigg)\bigg| \nonumber\\
&& \ \ \ \ = \ \bigg| \mathbf{P} \bigg( \frac{   \widetilde{H}_n \tau_n  }{ \sigma \sqrt{n} } \leq \frac{(u - \widetilde{K}_n)\tau_n}{\sigma \sqrt{n}a_n } \, \bigg| \,  \mathcal{F} \bigg)  -  \mathbf{P} \bigg( \mathcal{N}_1     \leq \frac{(u - \widetilde{K}_n)\tau_n}{\sigma \sqrt{n}a_n }  \, \bigg| \, \mathcal{F}  \bigg)\bigg| \ \leq  \frac{ C}{n^{\rho/2}} , \nonumber
 \end{eqnarray}
 where $C$ is a positive constant depending on $\sigma^2$ and  $\mathbf{E} Z_1^{2+ \rho}$.
Hence,  for   $ \rho \in (0, 1],$
\begin{eqnarray}\label{gcfsa1}
  P_1 \ \leq \frac{C }{n^{\rho/2}}.
\end{eqnarray}
On another hand, since $\mathcal{N}_1$ is independent of $ \widetilde{K}_n$, we have
\begin{equation}\label{gcfsa2}
P_2 =  \sup_{u \in \mathbf{R}} \bigg| \mathbf{P} \bigg(   \widetilde{K}_n  \leq u \bigg) -  \Phi \bigg(  \frac{ u  \tau_n }  { \sqrt{v_n}} \bigg) \bigg |
= \mathbf{K}\bigg(\frac{\widetilde{K}_n\tau_n }{\sqrt{v_n}}  \bigg) \ll  \frac{1}{ \sqrt{n}} + \frac{1}{v_n}  ,
 \end{equation}
where the last inequality comes from Theorem  \ref{thmBEERW}. Starting from (\ref{fsdop}) and considering the upper bounds \eqref{gcfsa1} and  \eqref{gcfsa2}, we obtain the desired estimate.

\section{Proof of Theorem  \ref{thm1ERW}}
\setcounter{equation}{0}
Let   $(N_i)_{1\leq i \leq n }$ be  a sequence of   $ \mathcal{N }(0,    a_i^2  )$-distributed independent random
variables. Assume moreover that  $(N_i)_{1\leq i \leq n }$ is independent of $(\Delta M_{i})_{1 \leq i \leq n }$ (recall \eqref{fsdfsdf} for the definition of $\Delta M_{i}$). For $n \geq3$, set $U_k=   \sum_{j=1}^{k}N_{j}, 1\leq k \leq n$.

As in the proof of Theorem \ref{thmBEERW}, we shall divide the proof in two cases according to the values of $p$.

\noindent \textit{1) Case  $p \in (0,3/4)$}. We shall use the same notations as those defined in  \eqref{defkappan} and  \eqref{deftkn}, selecting $\beta =1$. So
\begin{equation} \label{defkappatknW}
\kappa_n^2=\max_{[n/2]\leq i \leq n} a_i^2 \, , \,  t_{k,n}= \big ( \kappa_n^2  +  v_n-v_k \big )^{1/2}  \, , \,  u_{k,n}= 1  + 2 \sum_{i=k+1}^n  a_i^2 / \kappa_n^2  \, ,  \, 0 \leq k \leq n  .
\end{equation}
For any $r  >0$ and any $p \in [0,1]$,  let also
\begin{equation*} \label{defBnpF}
B_n (r,p)   =  \frac{1}{n^{r /2}}   +   \frac{1}{v_n }  .
\end{equation*}
Let  $N$ be a $ \mathcal{N}(0,  \kappa_n^2)$-distributed random variable independent of $(\Delta M_{i})_{1 \leq i \leq n }$  and $(N_i)_{1 \leq i\leq n }$.
Using Lemma 5.1 in Dedecker  \emph{et al.} \cite{DMR09}, for any $r \in (0, 2]$,  we first write
\begin{eqnarray}\label{gfgdfd01}
\zeta_r(\mathbf{P}_{M_n}, \mathbf{P}_{U_n}) \leq 2 \zeta_r(\mathbf{P}_{M_n+N}, \mathbf{P}_{U_n+N})  + 4\sqrt{2} \kappa_n^r .
\end{eqnarray}
Note that $\kappa_n^r \ll a_n^r \ll v_{n}^{r/2} /n^{r/2}$.
It remains to give an estimation for
\[
\zeta_r(\mathbf{P}_{M_n+N}, \mathbf{P}_{U_n+N}) = \sup_{f \in \Lambda_r}  \mathbf{E}[   f(M_n+N)-f(U_n+N)]   .
\]
With this aim, we shall use the Lindeberg method and denote
$$f_{n,k}(x)=\mathbf{E}\Big[   f\Big(x+U_n-U_k+N  \Big) \Big ]  \, ,  \, 1 \leq k \leq n  . $$
By independence of the sequences,
\begin{equation}  \label{LF1}
\mathbf{E}\Big[   f(  M_n +N)-f\Big(  \sum_{j=1}^{n}N_{j}+N\Big) \Big] = \sum_{k=1}^{n} D_k,
\end{equation}
where
$$D_k=  \mathbf{E}\Big[   f_{n,k}( M_{ k-1} + \Delta M_{ k})-f_{n,k}( M_{ k-1} + N_k) \Big] .    $$
Using twice a Taylor's expansion at order $5$,    we get
\begin{eqnarray} \label{LF2}
D_k =   I_{k,1}+I_{k,2}+I_{k,3}+I_{k,4}+I_{k,5},
\end{eqnarray}
where, for any integer $i \in [1, 4]$,
\[
 I_{k,i}= \frac{1}{i !} \mathbf{E}\Big[  f^{(i)}_{n,k}(M_{k-1}) \Delta M_{ k}^i    -f^{(i)}_{n,k}(M_{k-1})  N_k^i \Big] ,\\
\]
and
\[
   I_{k,5}= \frac1{120} \mathbf{E}\Big[  f^{(5)}_{n,k}(M_{k-1}+\vartheta_1 \Delta M_{k} ) \Delta M_{ k}^5    - f^{(5)}_{n,k}(M_{k-1} +\vartheta_2 N_k)  N_k^5 \Big] ,
\]
with $\vartheta_1, \vartheta_2$ some random variables with values in $ [0, 1].$
In the sequel, we give some estimations for $\sum_{k=1}^{n}I_{k,l}, l=1,2,\cdots, 5.$

Since $(M_k, \mathcal{F}_{k})$ is a martingale,   $I_{k,1}= \mathbf{E}\big[f'_{n,k}( M_{ k-1})\mathbf{E}[ \Delta M_{ k}- N_k |\mathcal{F}_{k-1}] \big ]=0.$
Therefore,\begin{eqnarray} \label{BI1nF}
\sum_{k=1}^{n}I_{k,1}=0.
\end{eqnarray}
We handle now the sum of  the $ | I_{k,2}| $'s. By Lemma 6.1 in \cite{DMR09}, for any positive integer $i$ such that $i \geq r$, there exists a positive positive constant $C_{r, i}$ such that for any integer $n\geq 1$ and any $k \in [1,n]$,
\begin{eqnarray}\label{devi}
  \| f^{(i)}_{n,k}(\cdot) \|_\infty \leq C_{r, i} t_{k,n}^{r-i}.
\end{eqnarray}
For $k=1$, we have
\begin{eqnarray}
   I_{1,2}  &=&   \mathbf{E}\Big[ \ \frac12 f''_{n,1}( M_{0})\Big( \mathbf{E}[\Delta M_{ 1}^2|\mathcal{F}_{0}]-  a_1^2   \Big) \Big] =0 . \nonumber
\end{eqnarray}
On another hand, by taking into account \eref{ddfsds}, (\ref{devi})  and \eqref{boundmoment2}, we get for $2\leq k \leq n,$
\begin{eqnarray}
  | I_{k,2} |  &=& \Big |   \mathbf{E}\Big[ \ \frac12 f''_{n,k}( M_{k-1})\Big( \mathbf{E}[\Delta M_{ k}^2|\mathcal{F}_{k-1}]-  a_k^2   \Big) \Big]  \Big |   \nonumber \\
  &\ll & ( 2p -1)^2 t_{k,n}^{r-2} \frac{a_k^2}{a_{k-1}^2} \mathbf{E}\Big[ \Big( \frac{ a_{k-1}T_{k-1}}{k-1} \Big)^2 \Big] \nonumber \\
  &\ll & ( 2p -1)^2 t_{k,n}^{r-2} \frac{a_k^2}{a_{k-1}^2}
   \frac{  (1+ v_{k-1}) }{ (k-1)^2 }. \label{B1Ik2F}
\end{eqnarray}
Clearly, when $p=1/2$, $  | I_{k,2}|  = 0$ for any positive integer $k$. In what follows we handle $\sum_{k=1}^n | I_{k,2} |$ when $p \in (0,3/4)\backslash \{1/2\}$. In this case,  by  (\ref{a10}) and (\ref{a15}),  we have $ \displaystyle \frac{  v_{k-1}}{ a_{k-1}^2(k-1)^2 } \asymp \frac{1  }{k}, k\rightarrow \infty$.
Therefore
\[
 \sum_{k=1}^{n} | I_{k,2} |   \ll  \max_{1\leq k \leq n/2}t_{k,n}^{r-2}
  \sum_{k \in [1, n/2]} \frac{ a_k^2 }{k} +  \frac{ 1}{n } \sum_{k\in (n/2, n]}  t_{k,n}^{r} \frac{a_k^2}{ t_{k,n}^2} .
\]
Clearly $ \max_{1\leq k \leq n/2}t_{k,n}^{r-2}    \leq C v_n^{(r-2)/2} $. Moreover, by \eqref{a10},
$$
\sum_{k \in [1, n]}  k^{-1} a_k^2   \ll  1 + (\log n) {\mathbf 1}_{p=1/2} + a^2_n  {\mathbf 1}_{p <1/2}   .
$$
 Hence, by taking in addition into account \eqref{ewer}, we derive  that, when $p \neq 1/2$,
\begin{equation*} \label{firstsecond}
 \sum_{k=1}^{n} | I_{k,2} | \ll v_n^{(r-2)/2}  ( 1+  a_n^2 {\bf 1}_{p < 1/2} )   +  \frac{ 1}{n } \sum_{k\in ( n/2, n]}  t_{k,n}^{r} \frac{a_k^2}{ t_{k,n}^2} .
\end{equation*}
Hence, taking into account \eqref{ine52cons} with $\alpha =r$, we derive that, for any $p$  in $(0,1/2) \cup (1/2,3/4)$,
\[
 \sum_{k=1}^{n} | I_{k,2} | \ll    v_n^{r/2}  ( v_n^{-1} + n^{-1})  +n^{-1}   \kappa_n^r  u^{r/2}_{[n/2],n} .
\]
But $ \kappa_n^r  u^{r/2}_{[n/2],n}  \ll v_n^{r/2}$.  So, overall,  when  $p \in (0,3/4)$ and $r \in (0,2]$,
\begin{equation} \label{I2BfinalF}
 \sum_{k=1}^{n} | I_{k,2} |  \ll  v_n^{r/2} B_n(r,p).
\end{equation}
We handle now the sum of the $|I_{k,3}| $'s.  By (\ref{devi}), Lemma \ref{lem2} and \eqref{boundmoment2}, we have
\begin{eqnarray*}
 | I_{k,3}|  &\leq& \frac16 \Big |  \mathbf{E}\Big[  f'''_{n,k}(M_{k-1}) \mathbf{E} [\Delta M_{ k}^3| \mathcal{F}_{k-1} ]   - f'''_{n,k}(M_{k-1}) \mathbf{E} [ N_k^3| \mathcal{F}_{k-1} ]  \Big |  \Big]\\
   &\ll&  t_{k,n}^{r-3} \mathbf{E} \big |\mathbf{E} [\Delta M_{ k}^3| \mathcal{F}_{k-1} ] \big |
   \ll  t_{k,n}^{r-3} a^3_k    \mathbf{E}\bigg|\frac{T_{k-1}}{k -1 }   \bigg| \\
  & \ll &   t_{k,n}^{r-3} \frac{a_k ^{3} }{(k-1)a_{k-1}  }  \sqrt{\mathbf{E}( a_{k-1} T_{k-1})^2}
\ll   t_{k,n}^{r-3}a_k ^{3} \frac{ \sqrt{v_{k-1}}}{(k-1)a_{k-1}  } ,
\end{eqnarray*}
with the convention $\displaystyle \sqrt{v_{0}}/ 0 a_{0}  =1.$
By  (\ref{a10}) and (\ref{a15}), and since  $p \in (0, 3/4)$, we have
$ \displaystyle  \frac{  \sqrt{v_{k-1}}}{ (k-1) a_{k-1} } =O\Big(\frac{ 1 }{\sqrt{k}} \Big),\ \  k \rightarrow \infty. $
Hence
\begin{equation} \label{decI3F}
  \sum_{k=1}^{n}   | I_{k,3} | \ll  A_n +  \frac{1}{\sqrt n}B_n  ,
\end{equation}
where
\[
A_n = \sum_{k \in [1, n/2)} \frac{ t_{k,n}^{r-3}a_k ^{3}  }{\sqrt{k}} \quad  \text{and} \quad B_n = \sum_{k \in ( n/2, n]} t_{k,n}^{r-3}a_k ^{3}  .
\]
By \eqref{a10}, it holds
\[
A_n  \ll  v_n^{(r-3)/2} \sum_{k  =1}^n \frac{a_k ^{3}  }{\sqrt{k}}  \ll   \left\{ \begin{array}{ll}
\displaystyle     v_n^{(r-3)/2}       \ \ \ & \textrm{if $ p>7/12 $}   \\
\nonumber \\
  v_n^{(r-3)/2}   \log n      \ \ \ & \textrm{if $ p=7/12$}   \\
\nonumber \\
    v_n^{(r-3)/2}  \sqrt{n} a_n^3     \ \ \ & \textrm{if $p < 7/12$}.\\
\end{array} \right.
\]
By \eqref{a10} and  \eqref{a15}, when $p\in (0,3/4)$, $v_n^{-3/2}  \sqrt{n} a_n^3  \ll  n^{-1}$ and
$
v_n^{-3/2} \log n  \ll   v_n^{-1} $.
Therefore, for any  $p \in (0, 3/4),$
\begin{equation} \label{BAnF}
A_n \ll   v_n^{r/2} B_n (r,p)  \, .
\end{equation}
We handle now the term $B_n$ in \eqref{decI3F}.  By taking into account \eqref{ine52cons}, we derive
\[
B_n \leq  \kappa_n {\mathcal T}_n (r-1)  \ll  \kappa^r_n  {\bf 1}_{(0,1) }(r) + \kappa_n ( \log n)  {\bf 1}_{r=1 } +  \kappa_n v_n^{(r-1)/2} {\bf 1}_{(1, 2] }(r)    \, .
\]
When, $p$ in $(0, 3/4),$ we have $\kappa^r_n \ll v_n^{r/2} n^{-r/2}$, $\kappa_n v_n^{(r-1)/2} \ll   v_n^{r/2} n^{-1/2} $ and $ \kappa_n ( \log n)
\ll   v_n^{1/2} n^{-1/2} ( \log n)$. Therefore, for any $p$ in $(0, 3/4),$
\begin{equation} \label{BBnF}
\frac{B_n}{\sqrt{n}} \ll   v_n^{r/2} B_n (r,p)  \, .
\end{equation}
So, starting from \eqref{decI3F} and taking into account \eqref{BAnF} and \eqref{BBnF}, we derive that, for any $r$ in $(0,2]$ and any $p$ in $(0, 3/4),$
\begin{equation} \label{BI3nF}   \sum_{k=1}^{n}   | I_{k,3} | \ll  v_n^{r/2} B_n (r,p)  \, .
\end{equation}

Before dealing with the sum of the $ | I_{k,4} | $'s, let us give an upper bound for  $ \sum_{k=1}^{n}   | I_{k,5} | $.  By (\ref{devi}) and Lemma \ref{lem1} again, we have
\[
|  I_{k,5} |  \ll  t_{k,n}^{r-5}\, \mathbf{E} \big[\mathbf{E} [ |\Delta M_{ k}|^5| \mathcal{F}_{k-1} ]+ \mathbf{E}  |N_k|^5   \big ]
 \ll t_{k,n}^{r-5}\, a^5_k  .
\]
Thus, by taking into account \eqref{ine52cons}, we have for $p \in (0, 3/4)$ and  $r \in (0,2]$,
\[
 \sum_{k=1}^{n}  | I_{k,5}  | \ll  v_n^{(r-5)/2}  \sum_{k=1 }^{[n/2]} a^5_k  +   \kappa^3_n \sum_{k=[n/2] +1 }^n  t_{k,n}^{r-5} a_k^2
  \ll   v_n^{(r-5)/2}  \sum_{k=1 }^{[n/2]} a^5_k + \kappa^r_n \, .
\]
Hence, by simple computations,
\begin{equation} \label{BI5nF}
 \sum_{k=1}^{n}  | I_{k,5}  |  \ll   v_n^{r/2} B_n (r,p)  \, .
\end{equation}
We handle now the term $ \sum_{k=1}^{n}   | I_{k,4} | $. By (\ref{devi}) and Lemma \ref{lem1}, we have
\[
|  I_{k,4} | = \frac1{24} \Big |  \mathbf{E}\Big[ f^{(4)}_{n,k}(M_{k-1}  ) \mathbf{E} [\Delta M_{ k}^4| \mathcal{F}_{k-1} ]   - f^{(4)}_{n,k}(M_{k-1} ) \mathbf{E} [ N_k^4]   \Big]  \Big | \ll  t_{k,n}^{r-4} a_k^{4} .
\]
Thus, by taking into account \eqref{ine52cons}, we have for $p \in (0, 3/4]$ and  $r \in (0,2]$,
\begin{multline*}
 \sum_{k=1}^{n}  | I_{k,4}  | \ll  v_n^{(r-4)/2}  \sum_{k=1 }^{[n/2]} a^4_k  +    \kappa^2_n \sum_{k=[n/2] +1}^n  t_{k,n}^{r-4}a_k^2 \\  \ll  v_n^{(r-4)/2}  \sum_{k=1 }^{[n/2]} a^4_k+ \kappa^r_n   \big (
  {\bf 1}_{r\in (0,2) } + ( \log n)    {\bf 1}_{r=2 }  \big ) \, .
\end{multline*}
 Thus, for $p \in (0, 3/4)$ and $r \in (0,2]$, we infer  that
 \begin{equation} \label{BI4nFinter}
 \sum_{k=1}^{n}  | I_{k,4}  |  \ll  v_n^{r/2} B_n (r,p)  {\bf 1 }_{r \neq 2} +   \kappa^2_n  ( \log n)  {\bf 1}_{r=2 }  .
\end{equation}
Therefore, starting from \eqref{LF1} and taking into account   \eqref{LF2}, \eqref{BI1nF} and   the upper bounds \eqref{I2BfinalF}, \eqref{BI3nF}, \eqref{BI5nF} and \eqref{BI4nFinter}, the result is proved for $r \in (0,2)$. It remains to prove it when $r=2$. With this aim, it is enough to prove that when $r=2$, $ \sum_{k=1}^{n}  | I_{k,4}  | \leq  C v_n B_n (2,p) $.
This means to delete the additional logarithmic term $\log n $ in the right-hand side of  \eqref{BI4nFinter}. From now on, we assume that $r=2$.  Note first that the previous computations proved that
$ \sum_{k=1}^{[n/2]}  | I_{k,4}  |  \ll  v_n B_n (2,p) $. So it remains only to prove that
\begin{equation} \label{BI4nFinterr2}
\sum_{k=[n/2]+1}^n  | I_{k,4}  |  \ll  v_n B_n (2,p) .
\end{equation}
With this aim, we first write the following decomposition: Since $ \mathbf{E} [ N_k^4 ]  = 3 a_k^4$,
\begin{eqnarray*}
|  I_{k,4} | &=& \frac1{24} \Big |  \mathbf{E}\Big[ f^{(4)}_{n,k}(M_{k-1}  ) \big (  \mathbf{E} [\Delta M_{ k}^4| \mathcal{F}_{k-1} ]   -  \mathbf{E} [ N_k^4 ]  \big ) \Big]  \Big |   \nonumber \\
&\leq &  I_{k,4}^{(1)} + I_{k,4}^{(2)} \, ,
\end{eqnarray*}
where
\[
 I_{k,4}^{(1)}:= \frac1{24}  \Big |  \mathbf{E}\Big[ f^{(4)}_{n,k}(M_{k-1}  ) \big (  \mathbf{E} [\Delta M_{ k}^4| \mathcal{F}_{k-1} ]   - a^4_k \big ) \Big]  \Big |
\text{ and }
 I_{k,4}^{(2)}:=   \frac{a_k^4}{12}   \Big |  \mathbf{E}\Big[ f^{(4)}_{n,k}(M_{k-1}  )  \Big]  \Big |   \, .
\]
Since $r=2$, by using  (\ref{devi}) and Lemma \ref{lem2}, we get that
\[
\sum_{k=[n/2]+1}^n I_{k,4}^{(1)}  \ll  \sum_{k=[n/2]+1}^n    t_{k,n}^{-2} a_k^{4} \frac{ \mathbf{E} (T^2_{k-1})}{(k -1)^2 }  \ll   \sum_{k=[n/2]+1}^n    t_{k,n}^{-2} \frac{a_k^{4} }{a_{k-1}^2} \frac{ v_{k-1}}{(k -1)^2 } .
\]
As $ k\rightarrow \infty$, recall that,  when $p \neq 3/4$, $ \displaystyle \frac{  v_{k-1}}{ a_{k-1}^2(k-1)^2 } \asymp \frac{1  }{k}$. Therefore
\[
\sum_{k=[n/2]+1}^n  I_{k,4}^{(1)}  \ll  \sum_{k=[n/2]+1}^n    t_{k,n}^{-2}  \frac{a_k^{4} }{k }   \ll    n^{-1}   \kappa_n^2  \sum_{k\in [ n/2, n]}  \frac{a_k^{2} } {t_{k,n}^2}    \, .
\]
Hence, taking into account \eqref{ine52cons}, we derive that
\[
\sum_{k=[n/2]+1}^n I_{k,4}^{(1)}  \ll   n^{-1}   \kappa_n^2 (\log n) \ll \kappa_n^2 . \]
 Note that when $p\in (0,3/4)$,  $  \kappa_n^2    \ll   v_n/ n \ll v_n B_n (2,p) $.
We then infer that  for $r=2$ and any $p \in (0,3/4)$,
\begin{equation*} \label{BI4nFinterr2B1}
\sum_{k=[n/2]+1}^n I_{k,4}^{(1)}  \ll   v_n B_n (2,p)    \, .
\end{equation*}
Therefore to end the proof of \eqref{BI4nFinterr2} and then of the theorem, we need to prove that
\begin{equation} \label{BI4nFinterr2B2}
\sum_{k=[n/2]+1}^n I_{k,4}^{(2)} \ll   v_n B_n (2,p)    \, .
\end{equation}
To achieve this, we need to have a better control of the term $\mathbf{E}\big[\big|f^{(4)}_{n,k}(M_{k-1}  )\big|\big]$ and we write
$$
\mathbf{E} \big[ f^{(4)}_{n,k}(M_{k-1}  )\big]= \mathbf{E} \big[ f^{(4)}_{n,k}(M_{k-1}  )\big]-\mathbf{E} \big[ f^{(4)}_{n,k}(G_{k-1})\big]+\mathbf{E} \big[ f^{(4)}_{n,k}(G_{k-1})\big]  ,
$$
where $G_{k-1}$ is a normal random variable with variance $a_1^2+ \cdots + a_{k-1}^2$.   By definition of $f_{n,k}$, we have
$$
\left |\mathbf{E}\big[  f^{(4)}_{n,k}(N_{k-1}) \big] \right |= \Big|  g_{n,k}^{(4)}(0) \Big|  \, ,
$$
where
\[
g_{n,k} (x) = { \mathbf E} \Big [ f \Big(x + \sum_{i=1, i \neq k}^n N_i + Z \Big) \Big] \, .
\]
Since, for any integer $k$  such that  $[n/2]+1 \leq k \leq n$, $ \kappa_n^2  + \sum_{i=1, i \neq k}^n a_i^2 \geq  v_n$, it follows from  (\ref{devi})  that $ |  g_{n,k}^{(4)}(0) | \ll  v^{-1}_n   $. Therefore,
$$
\left |\mathbf{E} \big[ f^{(4)}_{n,k}(G_{k-1})\big] \right | \ll  v^{-1}_n   .
$$
Hence, since $p \in (0,3/4)$,
\begin{equation} \label{ineq88}
\sum_{k=[n/2]+1}^n  a_k^4 \left |\mathbf{E} \big[ f^{(4)}_{n,k}(G_{k-1})\big]\right |   \ll   \frac{1}{v_n } \sum_{k=[n/2]+1}^n  a_k^4
\ll a_n^2
\ll n^{-1 }v_n .
\end{equation}
We now deal with the term
\begin{multline*}
\mathbf{E} \big[ f^{(4)}_{n,k}(M_{k-1}  )\big]-\mathbf{E} \big[ f^{(4)}_{n,k}(G_{k-1})\big]\ \ \ \ \\
=\int
\Big(\mathbf{E} \big[ f'(M_{k-1}  -u)\big] -\mathbf{E} \big[ f'(G_{k-1} -u)\big]\Big)\varphi_{t_{n,k}^2}^{(3)}(u) du. \ \ \ \ \  \ \ \ \
\end{multline*}
Hence, we have
\begin{multline*}
\left |\mathbf{E} \big[ f^{(4)}_{n,k}(M_{k-1}  )\big]-\mathbf{E} \big[ f^{(4)}_{n,k}(G_{k-1})\big] \right |\\
\leq c \, \sup_{u \in {\mathbf R}}\Big|\mathbf{E} \big[ f'(M_{k-1}  -u)\big] -\mathbf{E} \big[ f'(G_{k-1} -u)\big]\Big
|\int  \Big|\varphi_{t_{n,k}^2}^{(3)}(u) \Big| du. \ \ \ \ \ \ \
\end{multline*}
Recall that we have proved Theorem \ref{thm1ERW} for any $r \in (0,2)$. Hence applying it for $r=1$, we have, for any integer $k$ such that
$[n/2] +1 \leq k \leq n$,
$$
\sup_{u \in {\mathbf R}}\Big|\mathbf{E} \big[ f'(M_{k-1}  -u)\big] -\mathbf{E} \big[ f'(G_{k-1} -u)\big]\Big
 | \, \ll   \sqrt{\frac{v_n}{n} }  + \frac{1}{\sqrt{ v_n }}  .
$$
On another hand,
$$
\int  \Big|\varphi_{t_{n,k}^2}^{(3)}(u) \Big| du  \ll t_{k,n}^{-3}  .
$$
Consequently, by taking into account \eqref{ine52cons},
\begin{multline*}
\sum_{k=[n/2]+1}^n   a_k^4\Big|\mathbf{E} \big [ f^{(4)}_{n,k}(M_{k-1}  )\big]-\mathbf{E} \big [ f^{(4)}_{n,k}(G_{k-1})\big] \Big| \ll \kappa^2_n \Big (  \sqrt{\frac{v_n}{n} }  + \frac{1}{\sqrt{ v_n }} \Big )  \sum_{k=[n/2] +1}^n \frac{a_k^2}{t_{k,n}^3} \\ \ll  \kappa_n \Big (  \sqrt{\frac{v_n}{n} }  + \frac{1}{\sqrt{ v_n }} \Big )   \ll v_n   \Big (  \frac{1}{n}   + \frac{1}{ v_n \sqrt{ n }} \Big )   .
\end{multline*}
Hence, for any $p \in (0,3/4)$,
\begin{equation}
\sum_{k=[n/2]+1}^n   a_k^4\Big|\mathbf{E} \big [ f^{(4)}_{n,k}(M_{k-1}  )\big]-\mathbf{E} \big [ f^{(4)}_{n,k}(G_{k-1})\big] \Big| \ll   v_n  B_n(2,p)  . \label{ineq810}
\end{equation}

Combining the inequalities (\ref{ineq88}) and (\ref{ineq810}) together \eqref{BI4nFinterr2B2} follows. This completes the proof of the theorem when $p \in (0,3/4)$.

\smallskip

\noindent \textbf{2) Case $p=3/4$.} As in the proof of Theorem \ref{thmBEERW} in case $p=3/4$, we fix $\gamma \in (0,1)$ and we select $\kappa_n^2$ as in   \eqref{defkappanp=3/4}. The definitions of $ t_{k,n}$ and $ u_{k,n}$ are as in \eqref{defkappatknW}. Recall also that $v_n/\log n \rightarrow \pi/4$. We start again  from the inequality \eqref{gfgdfd01} and use the decompositions  \eqref{LF1} and \eqref{LF2}.  We still have to get upper bounds for the quantities $ | I_{k,i} |$ for $i =2, \cdots, 5$. For
$ | I_{k,2} |$ compared to the previous case, the difference is that $ \displaystyle \frac{  v_{k-1}}{ a_{k-1}^2(k-1)^2 } \asymp \frac{\log k  }{k}, k\rightarrow \infty$.  Hence starting from \eqref{B1Ik2F}, we get in this case that
\[ \sum_{k=1}^{n} | I_{k,2} |  \ll \max_{1\leq k \leq n^{\gamma} }t_{k,n}^{r-2}
  \sum_{k \in [1, n^{\gamma} ]} \frac{ \log k }{k^2} +  \frac{  \log n }{ n^{\gamma}  } \sum_{k\in ( n^\gamma , n]}  t_{k,n}^{r} \frac{a_k^2}{ t_{k,n}^2} \, .
\]
We then infer that when $r \in (0,2]$ and $p=3/4$,
\begin{equation*} \label{BoundI2p=3/4}
\sum_{k=1}^{n} | I_{k,2} | \ll  v_n^{(r-2)/2}  +  v_n^{r/2}   ( \log n)/  n^{\gamma}  \ll v_n^{(r-2)/2} .
\end{equation*}
We handle now the term $ \sum_{k=1}^{n}  | I_{k,3} | $.
When, $p =3/4$, the only difference  with the previous case is that
$ \displaystyle  \frac{  \sqrt{v_{k-1}}}{ (k-1) a_{k-1} } =O\Big(\frac{ \sqrt{ \log k } }{\sqrt{k}} \Big),\ \  k \rightarrow \infty$.  Hence,  \eqref{decI3F} has to be replaced by
\[
  \sum_{k=1}^{n}   | I_{k,3} | \ll \Big ( A'_n +  \frac{\sqrt{\log n }}{n^{\gamma/2}}B'_n \Big ) ,
\]
where  $ \displaystyle A'_n = \sum_{k \in [1, [n^\gamma]]} \frac{ t_{k,n}^{r-3}a_k ^{3}  \sqrt{ \log k} }{\sqrt{k}} $ and $B'_n = \sum_{k \in  ( n^\gamma, n] } t_{k,n}^{r-3}a_k ^{3} $.  Since $p=3/4$, $A_n' \ll v_n^{(r-3)/2}  $. On another hand, when $p=3/4$,  by using \eqref{ine52consp=3/4}, we infer that
\[
B'_n \ll \kappa_n^r  {\bf 1}_{(0,1) }(r) + \kappa_n (\log \log n)  {\bf 1}_{r=1 } + \kappa_n v_n^{(r-1)/2} {\bf 1}_{(1, 2] }(r) \, ,
\]
implying that $\frac{\sqrt{\log n }}{n^{\gamma/2}}B'_n \ll v_n^{(r-2)/2}$. Therefore, when $p=3/4$, for any $r$ in $(0,2]$,
\begin{equation*} \label{BoundI3p=3/4}  \sum_{k=1}^{n}   | I_{k,3} | \ll v_n^{(r-2)/2}   \, .
\end{equation*}
We handle now the term $ \sum_{k=1}^{n}   | I_{k,4} | $. We write this time
\[
\sum_{k=1}^{n}  | I_{k,4}  | \ll  v_n^{(r-4)/2} \sum_{k=1 }^{[n^{\gamma}]} a^4_k  +    \kappa^2_n \sum_{k=[n^{\gamma}] +1}^n  t_{k,n}^{r-4}a_k^2 \, .
\]
By using \eqref{ine52consp=3/4}, it follows that, for any $r \in (0,2]$,
 \begin{equation*}\label{BoundI3p=3/4}
 \sum_{k=1}^{n}  | I_{k,4}  | \ll  v_n^{(r-4)/2} + \kappa_n^r +  \kappa_n^2   (\log \log n ) {\mathbf 1}_{r= 2}  \ll v_n^{(r-2)/2}  \, .
\end{equation*}
Finally, let us handle  the term $\sum_{k=1}^{n}  | I_{k,5}  | $.  We write this time
\[
 \sum_{k=1}^{n}  | I_{k,5}  | \ll  v_n^{(r-5)/2}  \sum_{k=1 }^{[n^{\gamma}]} a^5_k  +   \kappa^3_n \sum_{k=[n^{\gamma}] +1 }^n  t_{k,n}^{r-5} a_k^2
.
\]
Thus, by taking into account \eqref{ine52consp=3/4}, we have for  $0< r \leq 2$,
\begin{equation*} \label{BoundI3p=3/4}
 \sum_{k=1}^{n}  | I_{k,5}  |  \ll   v_n^{(r-5)/2}  \sum_{k=1 }^{[n/2]} a^5_k + \kappa^r_n  \ll  v_n^{(r-2)/2}   \, .
\end{equation*}
Considering all these upper bounds, Theorem \ref{thmBEERW} is then proved in case $p=3/4$.

\section{Proof of Theorem  \ref{thm1}} \label{SectionPthm1}
\setcounter{equation}{0}
We start from the decomposition \eqref{dec1} of Section \ref{secp}. Let also
\begin{equation} \label{Notathm1F}
\sigma_{1,n}^2= \frac{n a_n^2 \sigma^2}{v_n + n a_n^2 \sigma^2}
\quad \text{and} \quad
\sigma_{2, n}^2 = \frac{v_n}{v_n + n a_n^2 \sigma^2} \, .
\end{equation}
 Let then $G_{1,n}$ and $G_{2, n}$ be two independent normal random variables with respective variances $\sigma_{1,n}^2$ and $\sigma_{2,n}^2$, and independent of $(\alpha_i, \beta_i, X_1, Z_i)_{i \geq 1}$.
By the triangle inequality,
\begin{equation}\label{b1bis}
\zeta_r(U_n+V_n, {\mathcal N}) \leq \zeta_r(U_n+V_n, G_{1,n}+V_n)+ \zeta_r(G_{1,n}+V_n, G_{1,n}+G_{2,n}) \, .
\end{equation}
We first deal with the first term on  right hand in \eqref{b1bis}. Let $\mathcal{F}=\sigma(\alpha_i, \beta_i, X_1, i\geq 1)$, and note that
\begin{equation}\label{b2bis}
\zeta_r(U_n+V_n, G_{1,n}+V_n)\leq {\mathbf E} [\zeta_r (P_{U_n+V_n|{\mathcal F}}, P_{G_{1,n}+V_n|{\mathcal F}})]\, .
\end{equation}
Now, since $V_n$ is ${\mathcal F}$-measurable and $G_{1,n}$ is independent of ${\mathcal F}$, we have
\begin{equation}\label{b3bis}
\zeta_r (P_{U_n+V_n|{\mathcal F}}, P_{G_{1,n}+V_n|{\mathcal F}})=\zeta_r (P_{U_n|{\mathcal F}}, P_{G_{1,n}}) \  \ a.s.
\end{equation}
Conditionally to ${\mathcal F}$, $U_n$ is a sum of independent random variables with variance
$a_n^2 \sigma^2/(v_n+ n a_n^2 \sigma^2)$ and absolute moment of order $2+\rho$ bounded by
$${\mathbf E}[\,|Z_1-1|^{2+\rho}\,] a_n^{2+\rho} /(v_n+ n a_n^2 \sigma^2)^{(2+\rho)/2}\, .$$ Applying Theorem 2.1 in \cite{DMR22} if $(r, \rho)\neq (1,1)$ and  the result of  Bikelis \cite{Bik66} if $(r, \rho)= (1,1)$, we infer that,  for any $r\in (0, 2]$ and $\rho \in (0, 1]$,
\begin{equation} \label{bikelisF}
 \zeta_r (P_{U_n|{\mathcal F}}, P_{G_{1,n}})\leq C  \left (\frac{\sqrt{n a_n^2 \sigma^2}}{\sqrt{v_n + n a_n^2 \sigma^2}}\right )^r
  \frac{1}{n^{(r\wedge\rho)/2}} \  \ a.s.
\end{equation}
where $C$ is a constant depending only on $r$, $\sigma^2$ and ${\mathbf E} Z_1 ^{2+\rho} $.
Hence, by  \eqref{b2bis},  \eqref{b3bis} and \eqref{bikelisF}, we get
\begin{equation}\label{b4bis}
 \zeta_r(U_n+V_n, G_{1,n}+V_n) \ll  \frac{1}{n^{(r\wedge\rho)/2}} .
\end{equation}

We now deal with the second term on the right hand side of \eqref{b1bis}. Since $G_{1,n}$ is independent of $(V_n, G_{2,n})$, it is easy to see that
\begin{equation*}\label{b5bis}
 \zeta_r(G_{1,n}+V_n, G_{1,n}+G_{2,n}) \leq  \zeta_r(V_n, G_{2,n}) \, .
\end{equation*}
Notice that
\begin{eqnarray*}
 \zeta_r(V_n, G_{2,n})&=& \left(\frac{ \sqrt{v_n}}{\sqrt{v_n +n a_n^2 \sigma^2}}\right)^r \zeta_r\!\left (\frac{a_n T_n-(2q-1)}{\sqrt {v_n}},
{\mathcal N} \right) \nonumber  \\
&\leq&   \  \zeta_r\!\left (\frac{a_n T_n-(2q-1)}{\sqrt {v_n}},
{\mathcal N} \right) . \label{b6bis}
\end{eqnarray*}
Thus from Theorem \ref{thm1ERW}, we get for $r \in (0, 2]$,
\begin{equation} \label{b7bis}
\zeta_r(G_{1,n}+V_n, G_{1,n}+G_{2,n})  \ll \frac{1}{n^{r/2}} + \frac{1}{v_n} .
\end{equation}
Starting from \eqref{b1bis} and considering the upper bounds  \eqref{b4bis}  and  \eqref{b7bis}, the result follows.

\section{Proof of Theorem  \ref{thm2}}
\setcounter{equation}{0}
We start again from the decomposition \eqref{dec1} of Section \ref{secp} and use the notations of Section \ref{SectionPthm1}.
By the triangle inequality
\begin{equation}\label{b1}
W_r(U_n+V_n, {\mathcal N}) \leq W_r(U_n+V_n, G_{1,n}+V_n)+ W_r(G_{1,n}+V_n, G_{1,n}+G_{2,n}) \, .
\end{equation}
We first deal with the first term on the right hand side of \eqref{b1}. By Fact 1.1 in \cite{DMR22sub}, note that
 \begin{equation}\label{b2}
W_r(U_n+V_n, G_{1,n}+V_n)\leq {\mathbf E} [W_r (P_{U_n+V_n|{\mathcal F}}, P_{G_{1,n}+V_n|{\mathcal F}})]\, .
\end{equation}
Now, since $V_n$ is ${\mathcal F}$-measurable and $G_{1,n}$ is independent of ${\mathcal F}$, we get
\begin{equation}\label{b3}
W_r (P_{U_n+V_n|{\mathcal F}}, P_{G_{1,n}+V_n|{\mathcal F}})=W_r (P_{U_n|{\mathcal F}}, P_{G_{1,n}}) \  \ a.s.
\end{equation}
Conditionally on ${\mathcal F}$, $U_n$ is a sum of independent random variables with variance
$a_n^2 \sigma^2/(v_n+ n a_n^2 \sigma^2)$ and absolute moment of order $2+\rho$ bounded by
$${\mathbf E}[\,|Z_1-1|^{2+\rho}\,]\, a_n^{2+\rho} /(v_n+ n a_n^2 \sigma^2)^{(2+\rho)/2}\, .$$ Applying Corollary 1.2 in \cite{Bob18}, we get  for any $r\in [1,2]$ and $\rho \in (0, r]$,
\begin{eqnarray}\label{b4}
 W_r (P_{U_n|{\mathcal F}}, P_{G_{1,n}})&\leq& C  \frac{\sqrt{n a_n^2 \sigma^2}}{\sqrt{v_n + n a_n^2 \sigma^2}}
 \left(\frac{{\mathbf E}[|Z_1-1|^{2+\rho}]}{\sigma^{2+\rho}}\right)^{1/r} \frac{1}{n^{\rho/2r}} \nonumber \\
 &\leq&\frac{C}{n^{\rho/2r}} \  \ a.s. \nonumber
\end{eqnarray}
where $C$ is a constant depending only on $r$, $\sigma^2$ and ${\mathbf E} Z_1 ^{2+\rho} $.
From (\ref{b2}) and \eqref{b3}, we get
\begin{equation} \label{b4.5}
W_r(U_n+V_n, G_{1,n}+V_n) \ll  \frac{1}{n^{\rho/2r}}\, .
\end{equation}

We now deal with the second term on the right hand side of \eqref{b1}. Since $G_{1,n}$ is independent of $(V_n, G_{2,n})$, it is easy to see that
\[
W_r(G_{1,n}+V_n, G_{1,n}+G_{2,n}) \leq W_r(V_n, G_{2,n}) \, .
\]
Now
$$
W_r(V_n, G_{2,n})= \frac{ \sqrt{v_n}}{\sqrt{v_n +n a_n^2 \sigma^2}}W_r\left (\frac{a_n T_n-(2q-1)}{\sqrt {v_n}},
{\mathcal N} \right) \, .
$$
By using inequality \eqref{Rioine},  we get that
$$
W_r(V_n, G_{2,n}) \leq c_r  \left (\zeta_r \left (\frac{a_n T_n-(2q-1)}{\sqrt {v_n}},
{\mathcal N} \right) \right )^{1/r}\, .
$$
Hence, from Theorem \ref{thm1ERW}, it follows  that, for any $r \in [1, 2],$
\begin{equation} \label{b6}
W_r(G_{1,n}+V_n, G_{1,n}+G_{2,n}) \ll \frac{1}{\sqrt{n}} + \frac{1}{v_n^{1/r}}  .
\end{equation}
The result follows by considering \eqref{b1} together with  \eqref{b4.5} and \eqref{b6}.


\end{document}